\documentclass[preprint,11pt]{amsart}%
\usepackage{graphicx}
\usepackage{amsmath,amssymb,amsfonts}
\usepackage{mathrsfs}
\usepackage[numbers]{natbib}
\usepackage[latin1]{inputenc}
\usepackage{graphicx,color}
\usepackage{tikz}
\usepackage{dutchcal}
\usepackage{caption}
\newtheorem{theorem}{Theorem}[section]

\newtheorem{lemma}[theorem]{Lemma}
\newtheorem{corollary}[theorem]{Corollary}

\theoremstyle{definition}
\newtheorem*{definition}{Definition}

\newtheorem{proposition}[theorem]{Proposition}
\theoremstyle{remark}
\newtheorem{remark}[theorem]{Remarks}

\numberwithin{equation}{section}

   \newcommand{\vstrut}{\rule{0in}{.16in}}
    
  \renewcommand{\Re}{\mbox{\rm Re}}
  \renewcommand{\Im}{\mbox{\rm Im}}
   \newcommand{\ph}{\mbox{$\varphi$}}

    \renewcommand{\phi}{\varphi}
   \newcommand{\Ht}{\mbox{$H^{2}$}}

   \newcommand{\D}{\mbox{$\mathbb{D}$}}
   \newcommand{\R}{\mathbb{R}}
  \newcommand{\Arg}{\text{Arg}}

    \newcommand{\BTC}{\widetilde{C_\phi}}
   \newcommand{\C}{\mbox{$C_{\varphi}$}}
   \newcommand{\CP}{\mathbb{C}}

\renewcommand{\phi}{\varphi}
\renewcommand{\Re}{{\rm Re\,}}
\renewcommand{\Im}{{\rm Im\,}}
\newcommand{\N}{\mathbb{N}}

\def\and{{\quad\text{and}\quad}}

\begin{document}

\title[Berezin Ranges of Composition Operators on $H^2$]{Berezin ranges of composition operators on $H^2$}
\author{Paul Bourdon}
\address{Department of Mathematics, University of Virginia,
         Charlottesville, USA}
\email{psb7p@virginia.edu}
\author{Trieu Le}
\address{Department of Mathematics and Statistics, University of Toledo, Toledo, OH, USA}
\email{Trieu.Le2@utoledo.edu}
\author{Derek Thompson}
\address{Department of Mathematics, Taylor University,  Upland, IN, USA}
\email{ theycallmedt@gmail.com}

%%% \date{August 24, 2022}
\date{October 8, 2026}

\begin{abstract} We explore properties of the Berezin range $B(C_\phi)$ of a composition operator $C_\phi$ on the Hardy space $H^2$, including constraints on its location in the complex plane $\CP$, topological properties, shape, symmetry, and convexity. For example, we show that  $B(C_\phi) \subseteq \{z\in \CP: \Re(z) > 0 \ \text{and}\ |z| \le 1 + |\phi(0)|\}$, that $B(C_\phi)$ is closed if and only if $C_\phi$ is the identity operator, and that $B(C_\phi)$ containa all its boundary points except for $0$ if $C_\phi$ is compact. Addressing questions raised by Cowen and Felder, we show that $C_\phi$ can be a (nondegenerate) disk and characterize when $B(C_\phi)$ is convex if $\phi$ is linear.   
\end{abstract}
\maketitle

\section{Introduction}
Let \Ht\ denote the classical \textit{Hardy space}, the Hilbert space of analytic functions $\displaystyle f(z) = \sum_{n=0}^{\infty}a_n z^n$ on the open unit disk \D\ such that $$\|f\|^{2}=\sum_{n=0}^{\infty}|a_n|^{2}<\infty.$$ The Hardy space \Ht\ is a \textit{reproducing kernel Hilbert space}, meaning that for each $\alpha\in \D$ there is a corresponding  function $K_{\alpha}\in H^2$ which reproduces a function's value at $\alpha$ through the inner product: $\langle f,K_\alpha\rangle=f(\alpha)$ for each $f\in \Ht$.  Here, $\langle \cdot, \cdot\rangle$ denotes the \Ht\ inner product given by 
$$
\langle f, g\rangle = \sum_{n=0}^\infty a_n\overline{b_n}, 
$$
where $(a_n)$ and $(b_n)$ are the Maclaurin-coefficient sequences of the $H^2$ functions $f$ and $g$, respectively. It is easy to check that $K_\alpha(z) = 1/(1-\bar{\alpha} z)$ and  $k_\alpha = \sqrt{1-|\alpha|^2}/(1-\bar{\alpha}z)$ is its normalization. 

 Every analytic selfmap $\phi$  of $\D$ induces a \textit{composition operator} \C\ on \Ht\ given by $\C f = f \circ \ph$, which is bounded and linear (see, e.g., \cite[Section 1.3]{shapirocomposition}). It is easy to verify that for a composition operator $C_\phi$ on \Ht\ and $\alpha\in \D$, $C_\phi^*K_\alpha = K_{\phi(\alpha)}$; in fact, a bounded linear operator $T$ on \Ht\ is a composition operator if and only if the set of reproducing kernels for \Ht\ is invariant under $T^*$ (\cite[Theorem 1.4]{cowen1995composition}). Composition operators were given two extensive treatments in the 1990s (see \cite{cowen1995composition} and \cite{shapirocomposition}), both illustrating the fascinating interplay between function theory and operator theory generated by the study of these operators.  Throughout this paper, $\phi$ denotes an analytic selfmap of $\D$  and $C_\phi$ represents the corresponding composition operator on $H^2$.

Composition operators on the Hardy space \Ht\ provide an important source of examples in operator theory.  For instance, a natural question raised by work of Berger and Coburn \cite{BergerCoburn} is whether a non-compact operator on a reproducing kernel Hilbert space may have a Berezin transform that vanishes on the unit circle.  Nordgren and Rosenthal \cite{NordgrenRosenthal} show that the answer to this question is yes by providing examples of  non-compact composition operators on \Ht\ whose Berezin transforms vanish on $\partial \D$ \cite[Theorem 2.3]{NordgrenRosenthal}.  More recently, Cowen and Felder \cite{cowen2022convexity} raise questions about the Berezin ranges of compositions operators,  which we address in this paper.   Our work also relates to that of Nordgren and Rosenthal's:  Theorem~\ref{ADT} below significantly generalizes \cite[Theorem 2.3]{NordgrenRosenthal}.

Recall that the numerical range of a bounded linear operator $T$ on a Hilbert space $\mathcal{H}$ is given by
$$ W(T) = \{ \langle Th, h \rangle : h \in \mathcal{H}, \|h\| =1\}.$$
If $\mathcal{H}$ is a reproducing kernel Hilbert space, then the {\it Berezin range} $B(T)$ of  $T: \mathcal{H} \to \mathcal{H}$ is the subset of $W(T)$ given by  $B(T) = \{ \langle Tk,k \rangle : k \in \mathcal{H}\ \text{is a normalized reproducing kernel}\}$.  Thus, the Berezin range of the bounded linear operator $T: H^2\to H^2$ is given by 
$$
B(T) = \{\langle T k_\alpha, k_\alpha \rangle: \alpha \in \D\}.
$$
The {\it Berezin transform} of $T:H^2\to H^2$ is the function $\widetilde{T}$ defined on $\D$ by
 $$
 \widetilde{T}(\alpha) = \langle T k_\alpha, k_\alpha \rangle;
 $$
 hence, the Berezin range of $T$ is the range of $\widetilde{T}: \D \to \CP$.  The {\it  Berezin radius} of $T$, denoted $b(T)$, given by $b(T) = \sup\{|z|: z\in B(T)\}$.

  In this paper, we explore  properties of Berezin ranges of composition operators on $H^2$, including constraints on  their location in the complex plane, topological properties, shape, symmetry, and convexity, with the latter three properties being tied to questions raised in \cite{cowen2022convexity}.
 For a composition operator $C_\phi:H^2\to H^2$, we have 
  $$
  \BTC(\alpha) = \langle C_\phi k_\alpha, k_\alpha \rangle= (1-|\alpha|^2)\left\langle \frac{1}{1-\bar{\alpha}\phi}, K_\alpha\right\rangle = \frac{1-|\alpha|^2}{1-\bar{\alpha}\phi(\alpha)}, \alpha\in \D.
  $$

     Because the linear-fractional transformation $z\mapsto 1/(1-z)$ maps the unit disk to the halfplane $\{z: \Re(z) > 1/2\}$, we see that for each $\alpha \in \D$, $\BTC(\alpha)$ is a point in this halfplane scaled by a positive factor $1-|\alpha|^2< 1$. Thus, the Berezin range $B(C_\phi)$ of a composition operator $C_\phi$ on \Ht\ is a subset of the open right halfplane $\{z: \Re(z) > 0\}$; in particular, 0 fails to be in $B(C_\phi)$ for all $\phi$.  Note also that $B(C_\phi)$ is both connected and bounded (by applying, respectively, the continuity of $\BTC$ and the Cauchy-Schwarz inequality, with the latter showing that $|z| \le \|C_\phi\|$ for $z\in B(C_\phi)$).

    Let $C_\phi$ be an arbitrary composition operator on $H^2$.    We begin the next section of our paper by showing that for any $\alpha\in \D$, we have $\left|\widetilde{C_\phi}(\alpha)\right| \le 1 + |\alpha||\phi(0)|$.  Thus,  $B(C_\phi)$ is contained in the half-disk $\{z: \Re(z) > 0\ \text{and}\ |z| \le 1 + |\phi(0)|\}$, and the Berezin radius $b(C_\phi)$ of $C_\phi$ is strictly less than $2$.

 If $\phi(z) = z$ is the identify function, then $C_\phi$ is the identity operator and, clearly, we have $B(C_\phi) = \{1\}$.     By the proof of  \cite[Theorem 3.1]{bourdon2002zero}, $0$ is in the closure of the Berezin range of every composition operator except for the identity operator. Because $0$ never belongs to the Berezin range of a composition operator, it follows that the only composition operator whose Berezin range is a closed subset of $\CP$ is the identity operator. In the next section of our paper, in addition to discussing Berezin radii of composition operators, we provide angular-derivative conditions on $\phi$ ensuring that $B(C_\phi)$ contains all of its boundary points except $0$.

 In Section 3, we show that if $\phi$ does not satisfy the angular-derivative conditions discussed in Section 2, then $B(C_\phi)$ might be an open set or might be a set containing some of its nonzero boundary points but not others.  We accomplish this by providing complete descriptions of $B(C_\phi)$ for $\phi$ belonging to the following two families of selfmaps of $\D$: 
$$
\mathcal{F}_1: = \{\phi: \phi(z) = a z + b\ \text{and}\ \phi\ \text{fixes a point on} \ \partial \D\},
$$
and 
$$
\mathcal{F}_2: = \left\{\phi_s: \phi_s(z) = \frac{z}{s -(s-1) z},\ 
\text{where}\  s >1\right\}.
$$
For $\phi\in \mathcal{F}_1$ for which $b\ne 0$, we show that $B(C_\phi)$ is the open disk with center $1/(1+a)$ and radius $1/(1+a)$, with the special case $\phi(z) =z/2+1/2$ answering a question raised by Cowen and Felder \cite[Question 5.4 \& Figure 5]{cowen2022convexity}. Note that if $\phi\in \mathcal{F}_1$ and $b= 0$, then $a$ is necessarily equal to $1$, so that $\phi(z) = z$ and $B(C_\phi) = \{1\}$. For each $\phi\in \mathcal{F}_2$, we show there is a corresponding number $\gamma\in (0, 1)$ such that for $z\in \partial B(C_\phi)$, we have $z\in B(C_\phi)$ if $\Re(z) > \gamma$ and $z\not\in B(C_\phi)$ if $\Re(z) \le \gamma$. We also provide in Section 3 computer-generated images of $B(C_\phi)$ for several $\phi$, focusing principally on selfmaps $\phi$ of the form $\phi(z) = az + b$.  These images suggest that when $\phi(z) = az +b$  does not fix a point in $\partial \D$, $B(C_\phi)$  can be convex or non-convex and that $0$ is a ``corner point'' of  $B(C_\phi)$, which we prove via a general sufficient condition  for $0$ to be a corner point of the Berezin range of a composition operator on $H^2$. We also show that $\phi(0) =0$ is both necessary and sufficient for $1$ to be a corner point of $B(C_\phi)$.

 In Section~\ref{Section4}, we describe $B(C_\phi)$   when $\phi(z) = az +b$.  In particular, we show that when $b\ne0$ and $\phi$ does not fix a point of $\partial \D$, then 
 $$
 B(C_\phi) = \left\{z\in \CP\setminus\{0\}: d \leq \Big|\frac{1}{z}-a\Big| - \Big|\frac{1}{z}-1\Big|
\right\},
$$
where $d = \sqrt{|1-a|^2 - |b|^2}$.
It follows that $\{1/z: z\in \partial B(C_\phi)\}$ is a branch of a hyperbola.  Then, we exploit this connection between $B(C_\phi)$ and hyperbolas to prove that $B(C_\phi)$ is convex if and only if $2d|\Im(a)|  \leq |b|(1-|a|^2)$, a condition that characterizes convexity of $B(C_\phi)$ for all selfmaps of $\D$ of the form $\phi(z) = az + b$. This result partially answers Question 5.2 of \cite{cowen2022convexity}: for a linear-fractional selfmap $\phi(z) = (az + b)/(cz+ d)$, what are necessary and sufficient conditions for $B(C_\phi)$ to be convex?

\section{The Berezin radius and topological properties of the Berezin range}

Because $1\in B(C_\phi)$ for every composition operator $C_\phi$ on $H^2$, the Berezin radius $b(C_\phi)$ of $C_\phi$ always satisfies $b(C_\phi)\ge 1$.  The next proposition gives an upper bound on $b(C_\phi)$, which is sharp when $\phi(0) = 0$.  Its proof relies on the generalized Schwarz Lemma for analytic selfmaps $\phi$ of $\D$: 
$$
|\phi(z)| \le \frac{|z| + |\phi(0)|}{1 + |z||\phi(0)|},
$$
for every $z\in\D$ (see, e.g., \cite[Corollary 2.40]{cowen1995composition}).
\begin{proposition}\label{brgub} Suppose that $\phi$ is an arbitrary analytic selfmap of $\D$; then for each $\alpha\in \D$, we have
$$
\left|\widetilde{C_\phi}(\alpha)\right| \le 1 + |\alpha| |\phi(0)|,
$$
so that $b(C_\phi)\le 1 + |\phi(0)|$.  Thus, $b(C_\phi)< 2$ for all analytic selfmaps $\phi$  of $\D$; and  if $\phi(0) = 0$, then $b(C_\phi) = 1$.
\end{proposition}
\begin{proof} Let $\alpha\in \D$ be arbitrary.  We have
\begin{align*}
\left|\widetilde{C_\phi}(\alpha)\right|& = \frac{1-|\alpha|^2}{|1-\bar{\alpha}\phi(\alpha)|}\\
& \le  \frac{1-|\alpha|^2}{1-|\alpha||\phi(\alpha)|}\\
&\le  \frac{1-|\alpha|^2}{1-|\alpha|\frac{|\alpha| + |\phi(0)|}{1 + |\alpha||\phi(0)|}} \\
& =  1 + |\alpha||\phi(0)|.
\end{align*}
Because $\alpha\in \D$, we see that $b(C_\phi) \le 1 + |\phi(0)|$. Because $\phi(0)\in \D$, we see that $b(C_\phi)< 2$. Also, if $\phi(0) =0$, then $b(C_\phi) \le 1$; and because $1\in B(C_\phi)$, we have $b(C_\phi) =1$. \end{proof}

 Hence, for every composition operator $C_\phi$ on $H^2$, we have $1 \le b(C_\phi) \le 1 + |\phi(0)|$, and the upper bound can be arbitrarily close to $2$.  The next proposition slightly improves this range of possible of values of $b(C_\phi)$ when $\phi(0)\ne 0$. 
  
\begin{proposition}\label{brgto} Suppose $\phi$ is an analytic selfmap of $\D$ such that $\phi(0)\ne 0$, then $b(C_\phi) >1$.\end{proposition}
\begin{proof} Let $\theta_0$ be an argument of $\phi(0)$.  Then, $\lim_{r\to 0^+} e^{-i\theta_0}\phi(re^{i\theta_0}) = |\phi(0)|$, and hence,
\begin{align*}
\lim_{r\to 0^+}\dfrac{\widetilde{C_\phi}(re^{i\theta_0}) -1}{r}
& = \lim_{r\to 0^+}\frac{1}{r}\left(\frac{1-r^2}{1-re^{-i\theta_0}\phi(re^{i\theta_0})}-1\right)\\
& = \lim_{r\to 0^+} \frac{e^{-i\theta_0}\phi(re^{i\theta_0}) - r}{1-re^{-i\theta_0}\phi(re^{i\theta_0})} = |\phi(0)| > 0.
\end{align*}
 Consequently, for $r$ sufficiently small and positive, \[\Re\left(\dfrac{\widetilde{C_\phi}(re^{i\theta_0}) -1}{r}\right) > 0, \]
 which implies $b(C_\phi)\ge \left|\widetilde{C_\phi}(re^{i\theta_0})\right| \ge \Re\left(\widetilde{C_\phi}(re^{i\theta_0})\right) >1$, as desired.
\end{proof}

Propositions~\ref{brgto} and \ref{brgub} yield the following.
\begin{corollary} Suppose that $\phi$ is an analytic selfmap of $\D$. Then $b(C_\phi)= 1$ if and only if $\phi(0) = 0$.
\end{corollary}

We show below that  angular derivatives of $\phi$ can play a role in determining $b(C_\phi)$ and that they also play a role in determining the topological properties of $B(C_\phi)$.     By definition, an analytic selfmap $\phi$ of $\D$ has an {\it angular derivative} at $\omega\in \partial \D$ provided  there is some $\eta\in \partial \D$ such that 
 \begin{equation}\label{DAngD}
 z\mapsto \frac{\eta - \phi(z)}{\omega - z}\  \text{has a finite nontangential limit at}\  \omega,
 \end{equation}
  in which case the limit is denoted $\phi'(\omega)$ and is called the angular derivative of $\phi$ at $\omega$. Note that if (\ref{DAngD}) holds, then $\phi$ necessarily has nontangential limit $\eta$ at $\omega$.  Also, note that if $\phi$ has an angular derivative at $\omega$ and $\eta$ is the nontangential limit of $\phi$ at $\omega$, then, in particular, we have
  $$
  \phi'(\zeta) = \lim_{r\to 1^-} \frac{\eta - \phi(r\omega)}{\omega-r\omega}.
  $$
    The Julia-Carath\'{e}odory Theorem (see, e.g., \cite[p.\ 57]{shapirocomposition}), provides useful necessary and sufficient conditions for $\phi$ to have an angular derivative at a point $\omega\in \partial\D$. For a point $\zeta\in\partial\D$ and a function $f$ defined on $\D$, we let $\angle\lim_{z\to \zeta} f(z)$ denote the nontangential limit of $f$ at $\zeta$.
\begin{theorem}[Julia-Carath\'{e}odory Theorem] \label{JCT} Suppose that $\phi$ is an analytic selfmap of $\D$ and that $\omega\in \partial \D$. The following are equivalent:
 \begin{itemize}
 \item[{\rm (JC 1)}]  $\liminf_{z\to \omega} \frac{1-|\phi(z)|}{1-|z|} = \delta < \infty$ $($where the $\liminf$ is taken as $z$ approaches $\omega$  unrestrictedly in $\D )$;
 \item[{\rm (JC 2)}] $\angle \lim_{z\to \omega} \frac{\eta - \phi(z)}{\omega - z}$ exits $($is finite$)$ for some $\eta \in \partial\D$ $($i.e., $\phi$ has an angular derivative at $\omega$$)$;
 \item[{\rm (JC 3)}] $\angle \lim_{z\to \omega} \phi'(z)$ exists, and $\angle \lim_{z\to \omega} \phi(z) = \eta\in \partial \D$.
 \end{itemize}
 Moreover,
 \begin{itemize}
 \item[(i)] $\delta > 0$ in  {\rm (JC 1)},
 \item[(ii)] the boundary points $\eta$ in {\rm (JC 2)} and {\rm (JC 3)} are the same, and
 \item[(iii)] the limit of {\rm (JC 2)} equals $\angle \lim_{z\to \omega} \phi'(z)$ from {\rm (JC 3)}, with both equal to $\omega \bar{\eta}\delta$; thus, the limit infimum of {\rm (JC 1)} is $|\phi'(\omega)|$ and $|\phi'(\omega)|=\delta >0$.
 \end{itemize}
 \end{theorem}

  We say that an analytic selfmap $\phi$ of $\D$ has $\zeta\in \partial \D$ as a fixed point and write $\phi(\zeta) = \zeta$ provided that $\lim_{r\to 1^-}\phi(r\zeta) = \zeta$.  A consequence of (iii) of the Julia-Carath\'{e}odory Theorem is that if $\phi$ has an angular derivative at a fixed point $\zeta\in \partial\D$, then $\phi'(\zeta) > 0.$  Every selfmap $\phi$ of $\D$ that is not an elliptic automorphism (an automorphism fixing a point in $\D)$ has an attractive fixed point in the closure $\D^-$ of $\D$. 
 
 \begin{theorem}[Denjoy-Wolff Theorem]   Suppose that $\phi$ is an analytic selfmap of $\D$ that is not an elliptic automorphism of $\D$. Then $\phi$ has a fixed point $w\in \D^-$ such that the sequence $\phi^{[n]}$ of iterates of $\phi$ converges to $w$ uniformly on compact subsets of $\D$.   
 \end{theorem}
 
For a proof of the Denjoy-Wolff Theorem, see, e.g., \cite[Theorem 2.51]{cowen1995composition}.
  The attractive fixed point $w$ of the Denjoy-Wolff Theorem is uniquely determined and called the {\it Denjoy-Wolff point} of $\phi$ and is also characterized as follows: $w$ is a point in $\D^-$  such $\phi(w) = w$ and $|\phi'(w)| < 1$ if $w\in \D$ and $ 0 < \phi'(w) < 1$ if $w\in \partial \D$.

  If $\phi$ has Denjoy-Wolff point $w\in \partial \D$ and $\phi'(w) < 1$, then $\phi$ is said to be of {\it hyperbolic type}.  For selfmaps of $\D$ of hyperbolic type, the next proposition improves the lower bound of the range $1 \le b(C_\phi) \le 1 + |\phi(0)|$; in addition, we will use the following proposition to show that the upper bound $1 + |\phi(0)|$ for $b(\C)$ can be sharp when $\phi(0)\ne 0$.   
  
\begin{proposition}\label{BFPE} Suppose that $\phi$ is an arbitrary analytic selfmap of $\D$ having angular derivative $\phi'(\zeta)$ at some $\zeta\in \partial\D$ that is fixed by $\phi$. Then,
$$
\lim_{r\to1^-} \widetilde{C_\phi}(r\zeta) = \frac{2}{1 + \phi'(\zeta)}.
$$
\end{proposition}
\begin{proof}  We have
\begin{align*}
\lim_{r\to 1^-}  \widetilde{C_\phi}(r\zeta) & 
= \lim_{r\to 1^-}  \frac{1-r^2}{1- r\bar{\zeta} \phi(r\zeta)}\\
&
= \lim_{r\to 1^-}  \frac{1-r^2}{1- r + r\bar{\zeta}\left(\zeta - \phi(r\zeta)\right)}\\
      & = \lim_{r\to 1^-}  \frac{1+r}{1 + r\frac{\left(\zeta - \phi(r\zeta)\right)}{\zeta(1 - r)}   } \quad (|\zeta| = 1\implies \bar{\zeta} = 1/\zeta)\\
      & = \frac{2}{1 + \phi'(\zeta)},
      \end{align*}
      as desired.
     \end{proof}
     
     A hyperbolic automorphism $\phi$ of $\D$ that is real on the real axis will necessarily fix $1$ and $-1$ and have the form $\phi^{}_t(z) = \frac{(t+1)z + t-1}{(t-1)z + t + 1}$ for some positive number $t \ne 1$. (To see that $\phi_t$ must have this form, observe that if $\phi$ is a  hyperbolic automorphism of $\D$ fixing $-1$ and $1$ and $T$ is the mapping $T(z) = (1+z)/(1-z)$ of $\D$ onto the open right halfplane $\mathbb{RHP}:=\{z: \Re(z) > 0\}$, then $z\mapsto T\left(\phi(T^{-1}(z))\right)$ is a hyperbolic automorphism of $\mathbb{RHP}$ fixing $0$ and $\infty$, and thus has the form $z\mapsto t z$ for some positive number $t$ not equal to $1$.)   Note that $\phi_t'(1) = 1/t$ while $\phi_t'(-1) = t$. Thus, $1$ is the Denjoy-Wolff point of $\phi^{}_t$ when $t > 1$ and $-1$ is the Denjoy-Wolff point if $0 < t < 1$.   If $t > 1$, then
     \begin{align*}
      \frac{2}{1 +\phi_t'(1)}  & \le  b(C_{\phi^{}_t}) \ \quad (\text{by Proposition~\ref{BFPE}})\\
                                   &  \le 1 + |\phi^{}_t(0)| \quad (\text{by Proposition~\ref{brgub}})\\
                                    & = \frac{2}{1+1/t} \\ 
                                      &= \frac{2}{1+\phi_t'(1)}.
      \end{align*}
 Hence, if $t > 1$, we have $b(C_{\phi^{}_t}) = 1 + |\phi^{}_t(0)|$. A similar argument shows the same is true if $0 < t < 1$.  Thus, we see that the upper bound $1 + |\phi(0)|$ for $b(C_\phi)$ can be sharp when $\phi(0)\ne 0$. 
 
   Recall from the Introduction that $B(C_\phi)$ is  closed only if $\phi(z) = z$ because otherwise $0$ is a limit point of $B(C_\phi)$ that does not belong to $B(C_\phi)$.  Observe that if $\phi(0) = 0$ and $\phi$ is not the identity function, then $B(C_\phi)$ is neither open nor closed because in this case $1 = b(C_\phi)$ is a boundary point of $B(C_\phi)$ contained in $B(C_\phi)$ while $0$ is a boundary point of $B(C_\phi)$ that is not contained in $B(C_\phi)$.  This raises the question of whether $B(C_\phi)$ can be an open set. In the next section, we show that if $\phi(z) = az + b$, then $C_\phi$ is open if and only if $b\ne0$ and $\phi$ fixes a point of $\partial\D$, in which case $B(C_\phi)$ is the open disk centered at $1/(1+a)$ having radius $1/(1+a)$. Theorem~\ref{ADT} below  shows that if $\phi$ is continuous on $\D^-$, then $\phi$ must fix a point in $\partial\D$  and have an angular derivative at that point in order for $B(C_\phi)$ to be open.  
   
   We note that angular derivatives play an important role in the analysis of compactness of composition operators on $H^2$.  If \C\ is compact, then $\phi$ fails to have an angular derivative at every point of $\partial\D$ \cite[Theorem 2.1]{ShapiroTaylor} and the converse is true if $\phi$ is univalent or boundedly valent \cite[Theorem 3.10]{maccluershapiro} (cf.\ \cite[Corollary 3.6]{ShapiroAnnals}).  
   
    For any compact operator $T$ on $H^2$, the Berezin transform $\widetilde{T}$ vanishes on $\partial\D$; i.e., $\widetilde{T}(\alpha)\to 0$ as $|\alpha| \to 1^-$. This holds because compact operators take weakly convergent sequences to norm convergent ones and the   normalized reproducing kernels $k_\alpha$ tend to 0 weakly in $H^2$ as $|\alpha| \to 1^-$ (see, e.g., \cite[Theorem 2.17]{cowen1995composition}). In some settings, the vanishing of the Berezin transform of an operator $T$ on $\partial\D$ implies $T$ is compact (see, e.g., \cite{AxlerZheng, MSWick, ZorboskaBerezin}).  However, as noted in the Introduction, Nordgren and Rosenthal \cite{NordgrenRosenthal} showed this is not the case for composition operators on $H^2$, providing a class of examples in their Theorem 2.3. For the examples of \cite[Theorem 2.3]{NordgrenRosenthal}, the composition-operator symbol $\phi$ is continuous on $\D^-$, fixes no point on $\partial\D$, and can be written as $\phi(z) = z\psi(z)$, where $\psi$ is univalent.  The following theorem provides a much larger class of examples.   Each symbol of \cite[Theorem 2.3]{NordgrenRosenthal} satisfies (ii) of the following theorem, which does not require that $\phi$ have a special form (only that $\phi$ be continuous on $\D^-$) and allows boundary fixed points as long as $\phi$ fails to have an angular derivative at these points. Moreover, while the symbols of compact composition operators satisfy (i) of the following theorem, (i) also provides examples of symbols $\phi$ not continuous on $\D^-$ such that $C_\phi$ is not compact and $\widetilde{C_\phi}$ vanishes on $\partial \D$; see  \cite[Section 10.2]{shapirocomposition}, where a description is provided of an infinite Blaschke product $\phi$ having no angular derivatives.

\pagebreak

\begin{theorem}\label{ADT}
    Let $\phi$ be an analytic selfmap of $\D$. If 
 \begin{itemize}
 \item[(i)] there is no point on $\partial \D$ at which $\phi$ has an angular derivative, or
 \item[(ii)]  $\phi$ extends to be continuous on the closure $\D^-$ of $\D$ and $\phi$ fails to have a fixed point on $\partial \D$ at which $\phi$ has an angular derivative, then
 \end{itemize}
 $\BTC$ vanishes on $\partial\D$ and $B(C_\phi)$ has the form $F\setminus\{0\}$, where $F$ is a compact subset (of the closed right halfplane). Thus, if (i) or (ii) holds, then $B(C_\phi)$ is not open. 
    \end{theorem} 
\begin{proof} We begin by showing that if $\BTC$ vanishes on $\partial \D$, then $B(C_\phi)$ has the form $F\setminus\{0\}$, where $F$ is  compact.  To simplify the notation, we let $f$ denote the Berezin transform of $C_\phi$, so that $f$ is the continuous function on $\D$ defined by 
\begin{equation}\label{dob}
f(\alpha) = \frac{1-|\alpha|^2}{1-\bar{\alpha} \phi(\alpha)}.
\end{equation}
Suppose that $f$ vanishes on $\partial\D$, i.e., $f(\alpha) \to 0$ as $|\alpha| \to 1^-$; then $f$ extends to be continuous on the closed unit disk $\D^-$ and the image of $\partial\D$ under the extended map $f_1$ is $\{0\}$. Because $\D^-$ is compact and $f_1$ is continuous, $f_1(\D^-)$ is a compact set $F$ containing $0$.  We have
$$
B(C_\phi) = f_1(\D) = \left(f_1(\D) \cup\{0\}\right)\setminus \{0\} = F \setminus\{0\}.
$$
Note that, in particular, the preceding description of $B(C_\phi)$ shows that $B(C_\phi)$ is not open, because $B(C_\phi)$ contains all its boundary points except $0$.  

We show that each of (i) and (ii) implies that $f$ vanishes on $\partial\D$ by proving the contrapositive implications. Suppose that there is an $\epsilon > 0$ and a sequence $(\alpha_n)$ in $\D$ with $\lim_{n\to \infty}|\alpha_n| = 1$ such that $|f(\alpha_n)|>\epsilon$ for all $n$. Note that there is some $\omega\in \partial \D$ and a subsequence $(\alpha_{j_n})$ of $(\alpha_n)$ such that $\lim_{n\to\infty} \alpha_{j_n} = \omega$.  We have for all $n$,
\begin{align*}
    \frac{1-|\alpha_{j_n}|}{1-|\phi(\alpha_{j_n})|} & = \frac{1}{1+|\alpha_{j_n}|} \frac{1-|\alpha_{j_n}|^2}{|1-\overline{\alpha_{j_n}}\phi(\alpha_{j_n})|}\\
    &\ge \frac{1}{2} \frac{1-|\alpha_{j_n}|^2}{|1-\overline{\alpha_{j_n}}\phi(\alpha_{j_n})|}\\
    & = \frac{1}{2} |f(\alpha_{j_n})|\\
    & > \epsilon/2.
\end{align*}
Thus,
$$
\liminf_{z\to \omega} \frac{1-|\phi(z)|}{1-|z|} < 2/\epsilon,
$$
so that $\phi$ has an angular derivative at $\omega$ by Julia-Carath\'{e}odory Theorem. Thus, (i) implies $f$ vanishes on $\partial\D$.

  From the preceding paragraph, we know that $\phi$ has an angular derivative at $\omega\in \partial \D$ and that there is a sequence $(\alpha_{j_n})$  in $\D$ converging to $\omega$ such that the sequence
\begin{equation}\label{sbao} \left( \frac{1-|\alpha_{j_n}|^2}{|1-\overline{\alpha_{j_n}}\phi(\alpha_{j_n})|}\right)
\end{equation}
is bounded away from $0$.  Because $\phi$ has angular derivative at $\omega$, there is an $\eta\in \partial\D$ such that $\phi$ has radial limit $\eta$ at $\omega$ by the Julia-Carath\'{e}odory Theorem.  Now suppose that $\phi$ is continuous on $\D^-$ (thus, in particular, at $\omega$); then $\lim_{n\to\infty} \phi(\alpha_{j_n}) = \eta$ and because the sequence (\ref{sbao}) is bounded away from $0$ we see that $\eta = \omega$.  Thus, $\phi$  has radial limit $\omega$ at $\omega$; that is, $\phi$ has an angular derivative at a fixed point $\omega\in \partial \D$. Hence, (ii) implies $f$ vanishes on $\partial\D$.
\end{proof}

The necessary conditions of Theorem~\ref{ADT} for $B(\C)$ to be open ($\phi$ must have an angular derivative, or, if continuous on $\D^-$, must have an angular derivative at a fixed point in $\partial\D$) are not sufficient.  For example, $\phi(z) = z/(2-z)$ is an analytic selfmap of $\D$ that is continuous on $\D^-$ and has angular derivative 2 at its fixed point 1. However,  $B(\C)$ is not open because $\phi(0) = 0$, so that $1$ is a boundary point of $B(\C)$ that is contained in $B(\C)$.  In the next section, we completely describe $B(\C)$ when 
$\phi$ has the form  $\phi_s(z) = \frac{z}{s - (s-1)z}$, where $s>1$, showing that, e.g.,  $B(C_{\phi_2})$  contains precisely those boundary points whose real parts exceed 9/16.

We conclude this section with the following proposition whose proof relies on Theorem~\ref{ADT}, and Proposition \ref{brgto}. 

\begin{proposition} Suppose that $\phi$ is an analytic selfmap of $\D$ having Denjoy-Wolff point $w\in \D$ such that  $\phi$ is continuous on $\D^-$. Then the Berezin radius of $C_\phi$ is attained: there is a number $\alpha \in \D$ such that $\left|\widetilde{C_\phi}(\alpha)\right| = b(C_\phi)$. 
\end{proposition}
\begin{proof}  If the Denjoy-Wolff point $w$ is $0$, then $\phi(0) = 0$ and $b(C_\phi) = 1$ by Proposition~\ref{brgub}; thus,  $\widetilde{C_\phi}(0) = b(C_\phi)$, so that the Berezin radius of $C_\phi$ is attained.   
 
Let $w\in \D\setminus\{0\}$. Because the iterate sequence $(\phi^{[n]}(0))$ converges to $w$, we see that $\phi$ cannot fix $0$.  Thus, by Proposition~\ref{brgto}, $b(C_\phi) > 1$.  

Suppose, in order to obtain a contradiction that the Berezin radius of $C_\phi$ is not attained at a point in $\D$.  
 Then, there is a sequence $(\alpha_j)$ in $\D$ with $(|\alpha_j|)$ converging to $1$ such that 
\begin{equation}\label{btcbr}
\lim_{j\to\infty} \left|\widetilde{C_\phi}(\alpha_j)\right|  = b(C_\phi) > 1.
\end{equation}
Without loss of generality, we may assume that $(\alpha_j)$ converges to a point $\zeta\in \partial\D$. Note that $(\phi(\alpha_j))$ must also converge to $\zeta$ because otherwise $(\phi(\alpha_j))$ would have a subsequence $(\phi(\alpha_{k_j}))$ converging to a point of $\D^-$ different from $\zeta$ and the sequence $\left( \left|\widetilde{C_\phi}(\alpha_{k_j})\right|\right)$ would converge to $0$, contradicting (\ref{btcbr}). We conclude that $\phi(\zeta) = \zeta$ because $\phi$ is continuous on $\D^-$.  

From the definition of $\widetilde{C_\phi}$ and (\ref{btcbr}), we obtain
$$
\limsup_{j\to\infty} \frac{1-|\alpha_j|^2}{1-|\alpha_j||\phi(\alpha_j)|} \ge \lim_{j\to \infty}\left|\widetilde{C_\phi}(\alpha_j)\right| > 1.
$$
Thus,
\begin{align*}
1 & > \liminf_{j\to \infty} \frac{1-|\alpha_j||\phi(\alpha_j)|}{1-|\alpha_j|^2}\\ 
&= \liminf_{j\to \infty} \frac{1-|\alpha_j| +|\alpha_j|(1-|\phi(\alpha_j)|)}{1-|\alpha_j|^2} \\
& = \liminf_{j\to \infty} \frac{1 +|\alpha_j|\left(\frac{1-|\phi(\alpha_j)|}{1-|\alpha_j|}\right)}{1+|\alpha_j|},
\end{align*}
and it follows that 
$$
\liminf_{j\to\infty}\frac{1-|\phi(\alpha_j)|}{1-|\alpha_j|} < 1.
$$
Because $\lim_{j\to\infty} \alpha_j = \zeta\in \partial \D$ and $\phi(\zeta) = \zeta$, the Julia-Carath\'{e}odory theorem tells us that $\phi$ has an angular derivative at $\zeta$ and $0 < \phi'(\zeta) < 1$.  We conclude that $\zeta$ is the Denjoy-Wolff point of $\phi$; hence, we must have $\zeta = w$, a point on $\partial \D$ equalling a point in $\D$, a contradiction.  Thus, the Berezin radius of $C_\phi$ is attained.
\end{proof}

\section{Examples and additional properties of $B(\C)$ }\label{ExAddProp}

In this section, we provide complete descriptions of $B(C_\phi)$ for $\phi$ belonging to the following two families of selfmaps of $\D$: 
\begin{equation}\label{FamF1}
\mathcal{F}_1: = \{\phi: \phi(z) = a z + b\ \text{and}\ \phi\ \text{fixes a point on} \ \partial \D\}
\end{equation}
and 
\begin{equation}\label{FamF2}
\mathcal{F}_2: = \left\{\phi_s: \phi_s(z) = \frac{z}{s - (s-1)z},\ \text{where}\ s >1\right\}.
\end{equation}
We also provide a simple sufficient condition for the Berezin range of a composition operator to be symmetric about the real axis (addressing Question 5.3 of \cite{cowen2022convexity}) and introduce the notion of  ``corner point" for $B(C_\phi)$, providing a sufficient condition for $0$ to be a corner point as well as a necessary and sufficient condition for $1$ to be a corner point.

\subsection{$B(C_\phi)$ for $\phi\in \mathcal{F}_1$}
We begin with an elementary proposition describing members of the family $\mathcal{F}_1$ in terms of their coefficients $a$ and $b$. 

\begin{proposition}\label{LSMBFP} The selfmap $\phi(z) = az + b$ fixes a point on $\partial\D$ if and only if 
\begin{equation}\label{CFPB}
0 < a  \le 1\ \text{and}\  a + |b| = 1.
\end{equation}
\end{proposition}
\begin{proof} Let  $\phi(z) = az + b$ be a selfmap of $\D$.  Observe that $\phi$ cannot fix point on $\partial \D$ if $a=0$. Thus, we assume that $a\ne 0$.       Because $\phi$ is continuous on $\D^-$  and  a selfmap of $\D$ ,
\begin{equation}\label{CSF}
1 \ge \max \{|\phi(\zeta)|: \zeta\in\partial\D\}  = |a| + |b|,
\end{equation}
with the maximum attained at $\zeta = \bar{a}/|a|$ if $b = 0$ and at $\zeta = \bar{a}/|a|\cdot b/|b|$ if $b\ne 0$.
(Remark:  given our assumption that $a\ne 0$, it's easy to check that $|a| + |b| \le 1$ is sufficient for $\phi(z) = a z + b$ to be a selfmap of $\D$.)

Suppose that (\ref{CFPB}) holds.   If $a = 1$, then $\phi(z) = z$ and $\phi$ fixes all points of $\partial\D$. If $0< a < 1$, then $\phi$ fixes $b/(1-a)$, which is a point on $\partial \D$ because $|b/(1-a)| = |b|/(1-a)$ and $|b| = 1-a$ by (\ref{CFPB}).   

Conversely, suppose $\phi$ fixes a point $\zeta\in \partial\D$, so that $a\zeta + b = \zeta$. If $b = 0$, then $a = 1$ and condition (\ref{CFPB}) is satisfied.  Let $b \ne0$; then $a\ne 1$, and $\zeta = b/(1-a).$  Because $\zeta\in \partial \D$, we have $ |b| = |1-a|$. By (\ref{CSF}), we have $|b|  \le 1-|a|$, so that
$$
1-|a| \ge |b| = |1-a| \ge 1 - |a|,
$$
from which it follows that $|b| = 1- |a|$ and  $|1-a| = 1-|a|$, with the latter equality implying $0< a < 1$ because we are assuming $a\ne 0$ and $a\ne 1$.  Thus, we have shown that if $b\ne0$, then $a$ and $b$ satisfy (\ref{CFPB}), completing the proof.
\end{proof}

Our first theorem of Section~\ref{ExAddProp}, Theorem~\ref{diskunion} below, shows that for a linear selfmap $\phi(z) = a z + b$  of $\D$  with $b \ne 0$, the Berezin transform $\widetilde{C_\phi}$ maps each of the circles $\{z: |z| = r\}$, $0 < r < 1$, to a circle and that $B(C_\phi)$ is the union of these image circles and their interiors.   

\begin{lemma}\label{TLL}
Let $I$ be either $[\alpha,\beta]$ or $[\alpha,\beta)$ for some real numbers $\alpha<\beta$. Let $c:I\rightarrow\CP$ and $\rho: I\rightarrow [0,\infty)$ be continuous functions such that $\rho(\alpha)= 0$. Then
\[\bigcup_{r\in I}\Big\{z: |z-c(r)|\leq \rho(r)\Big\} = \bigcup_{r\in I}\Big\{z: |z - c(r)| = \rho(r)\Big\}.\]
\end{lemma}

\begin{proof}
It suffices to prove that the union on the left side is contained in union on the right side. Let $r_0\in I$ and let $z_0$ satisfy $|z_0 - c(r_0)| \leq \rho(r_0)$. Define the continuous real-valued function $g: I\rightarrow\R$ by $g(r) = |z_0 - c(r)| - \rho(r)$ for $r\in I$. Note that $g(\alpha) = |z_0-c(\alpha)|\geq 0$ and $g(r_0) = |z_0 - c(r_0)| - \rho(r_0) \leq 0$. By the Intermediate-Value Theorem, there exists $r_1\in [\alpha,r_0]$ such that $g(r_1) = 0$, which means that $z_0$ belongs to the union on the right side. This completes the proof of the lemma.
\end{proof}

\begin{theorem}\label{diskunion}
Suppose that $\phi(z) = az + b$ is a selfmap of $\D$ with $b\ne 0$; then the Berezin range of $C_\phi: H^2 \to H^2$ is a union of closed disks having positive radius, where the boundaries of these disks are precisely the images of the   circles $\{z: |z| = r\}$, $0 < r < 1$,  under $\BTC$.
\end{theorem}
    \begin{proof} Let $\phi(z) = az + b$ be a selfmap of $\D$ with $b\ne 0$. (Recall from the proof of Proposition~\ref{LSMBFP} that if $a\ne 0$, then $\phi(z) = a z + b$ is a selfmap of $\D$ if and only if $|a| + |b| \le1$.) We have 
\begin{align}
B(C_\phi)  & = \{\widetilde{C_\phi}(\alpha) : \alpha \in \D\} \nonumber\\ 
&= \left\{(1-|\alpha|^2)  \frac{1}{1-\bar{\alpha} \phi(\alpha)}: \alpha \in \D\right\} \nonumber\\ 
&= \left\{(1-|\alpha|^2)  \frac{1}{1 - a|\alpha|^2-b\bar{\alpha}}: \alpha \in \D\right\} \nonumber\\ 
&= \left\{(1-r^2)  \frac{1}{1 - ar^2-bre^{-i\theta}}: 0\le r < 1,  0\le \theta < 2\pi \right\}.\label{modbonly}
\end{align}
Observe that the final description of $B(C_\phi)$ can be viewed as the union of the circles  $$f_r\left(\vstrut \{z: |z| = 1\}\right), 0 \le r < 1,$$ where $f_r$ is the linear-fractional transformation $f_r(z) = \frac{1-r^2}{1-ar^2 - brz}$.  Note that for $0 < r < 1$, these circles have positive radius and are the images of the circles $\{z: |z| = r\}$ under the Berezin transform $\widetilde{C_\phi}$.

  For $0 \le r < 1$, let $c(r)$ be the center of the circle $f_r(\{z: |z| = 1\})$ and let $\rho(r)$  be the corresponding radius. It is easy to check that $c(0)=1$, $\rho(0) = 0$, $\rho(r) > 0$ for $r\in (0,1)$, and that $c$ and $\rho$ are continuous on $[0, 1)$. (For instance, $c(r)$ is the average of the diametrically opposed points 
$$
\frac{1-r^2}{1-ar^2}\cdot\frac{1}{1+\frac{|b|r}{|1-ar^2|}}\ \text{and}\  \frac{1-r^2}{1-ar^2}\cdot\frac{1}{1-\frac{|b|r}{|1-ar^2|}}
$$
on $f_r(\{z:|z| = 1\}$ and $\rho(r)$ is half the distance between them.)  Note our work of the preceding paragraph shows that
$$B(C_\phi)= \bigcup_{r\in[0,1)}    \{z: |z-c(r)| = \rho(r)\}.$$

For $r\in [0,1)$, let $F_r =\{z: |z-c(r)|\le \rho(r)\}$, so that for $r\in (0,1)$, $F_r$ is a closed disk of positive radius whose boundary is $\BTC(\{z:|z| = r\})=f_r(\{z: |z| =1\})$. We claim $1\in  \cup_{r\in (0,1)} F_r$.  Assuming this claim to be valid, we have
$$
B(C_\phi)= \bigcup_{r\in[0,1)}    \{z: |z-c(r)| = \rho(r)\} =   \left(\bigcup_{r\in [0,1)} F_r\right) = \{1\} \cup\bigcup_{r\in (0,1)} F_r = \bigcup_{r\in (0,1)} F_r,
$$
where the second equality follows from  Lemma~\ref{TLL}. Thus, our proof is complete once we have established the validity of our claim. For $r\in (0,1)$, we have  
\begin{align*}
\left|\frac{1-ar^2}{1-r^2}  - 1 \right| &= \frac{r^2|1-a|}{1-r^2} \\
&=\frac{r^2|1-a|}{1-r^2} -\frac{|b|r}{1-r^2} + \frac{|b|r}{1-r^2} \\ 
&= -r\frac{|b|-r|1-a|}{1-r^2} + \frac{|b|r}{1-r^2}\\
& < \frac{|b|r}{1-r^2}\quad \text{for}\ r\ \text{sufficiently small}.
\end{align*}
Hence, for $r_0>0$ sufficiently small, $1$ lies in the interior of the circle with center $\frac{1-ar_0^2}{1-r_0^2}$ and radius $\frac{|b|r_0}{1-r_0^2}$; that is, $1$ lies in the interior of the circle $\mathfrak{C}_{r_0}:=\left\{\frac{1-ar_0^2 - br_0e^{it}}{1-r_0^2}: 0 \le t < 2\pi\right\}$. Observe that 0 lies outside of $\mathfrak{C}_{r_0}$.   Thus,  $z\mapsto 1/z$  maps $\mathfrak{C}_{r_0}$ to  a circle containing $1$ in its interior. This completes the proof because the image of $\mathfrak{C}_{r_0}$ under $z\mapsto 1/z$ is $f_{r_0}(\{z: |z| = 1\})$, so that $1\in F_{r_0}$.
\end{proof}

\begin{corollary}\label{ab1}
    Suppose that $\ph(z)=az+b$  is a selfmap of $\D$, that $b\ne 0$, and that $\phi$ has a fixed point on $\partial \D$. Then $B(\C)$ is an open disk  with center $1/(1+a)$ and radius $1/(1+a$).  
\end{corollary}
\begin{proof}
  By Proposition~\ref{LSMBFP} and its proof, $0 < a  < 1\ \text{and}\  a + |b| = 1$.  By Theorem \ref{diskunion} and its proof, we see that $B(\C)$ is the union of closed disks $F_r$ with diameters extending from  $\alpha_r:=\frac{1-r^2}{1-ar^2 +|b|r} $ to $\beta_r:=\frac{1-r^2}{1-ar^2 -|b|r}$, where $0< r < 1$. It is easy to check that as $r$ increases to $1$,  $\alpha_r$ decreases to $0$ while $\beta_r$ increases to $2/(1+a)$ (keep in mind that $a + |b| =1$).  Thus, these closed disks are nested with $F_{r_1}\subset F_{r_2}$ when $0< r_1 < r_2 < 1$, and $B(C_\phi) = \cup_{r\in (0,1)} F_r$ is the open disk with center $1/(1+a)$ and radius $1/(1+a)$. 
\end{proof}

As noted in the Introduction, the preceding corollary gives a positive answer to a question in \cite{cowen2022convexity} by establishing that $B(\C)$ is a disk (more precisely, an open disk) when $\ph(z)=z/2+1/2$.  Corollary~\ref{ab1} also shows that Berezin ranges of composition operators on $H^2$ can be open sets.   Finally, note that any member of the family $\mathcal{F}_1$ (see (\ref{FamF1}) above) either satisfies the hypotheses of Corollary~\ref{ab1} or is the identity function $\phi(z) = z$, which induces a composition operator (the identity operator on $H^2$) whose Berezin range is $\{1\}$. Thus, as promised, we have provided a complete description of Berezin ranges of composition operators whose symbols belong to $\mathcal{F}_1$.

Now observe that for every $\phi\in \mathcal{F}_1$, $B(\C)$ is convex, being a disk (degenerate if $b =0$) in the complex plane. 
As noted in the Introduction, we prove in Section~\ref{Section4} that the selfmap $\phi(z) = az + b$ of $\D$ induces a composition operator $C_\phi$ on $H^2$ such that $B(C_\phi)$ is convex if and only if 
\begin{equation}\label{ConvexityCondition}
2|\Im(a)|d \le |b|(1-|a|^2), \ \text{where} \ d = \sqrt{|1-a|^2 - |b|^2}. 
\end{equation}
  Note that the preceding condition holds for $\phi\in \mathcal{F}_1$ because if $\phi(z) = az + b$ is in $\mathcal{F}_1$, then $\Im(a) = 0$. 
  
\begin{figure}[h]
    \centering
    \includegraphics[height=1.6in]{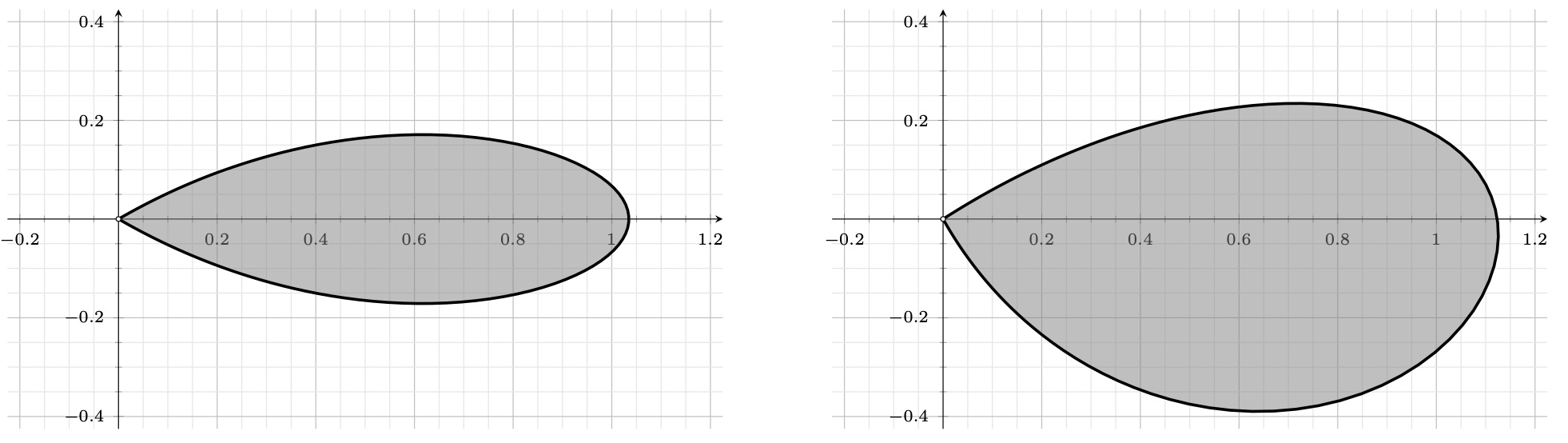}
    \caption{\small $B(C_\phi)$ for (i) $\phi(z) = \frac{1}{2}z + \frac{1}{4}$ (left) and (ii) $\phi(z) = \left(\frac{1}{3} - \frac{1}{6}i\right)z + \frac{1}{2}$. Both ranges are convex with the convexity condition~(\ref{ConvexityCondition}) satisfied: (i) $0 \le  \frac{3}{16}$; (ii) $ \frac{\sqrt{2}}{9} \le \frac{31}{72}$.}\label{ConFigs1}

\end{figure}

\begin{figure}[h]
    \centering
    \includegraphics[height=1.6in]{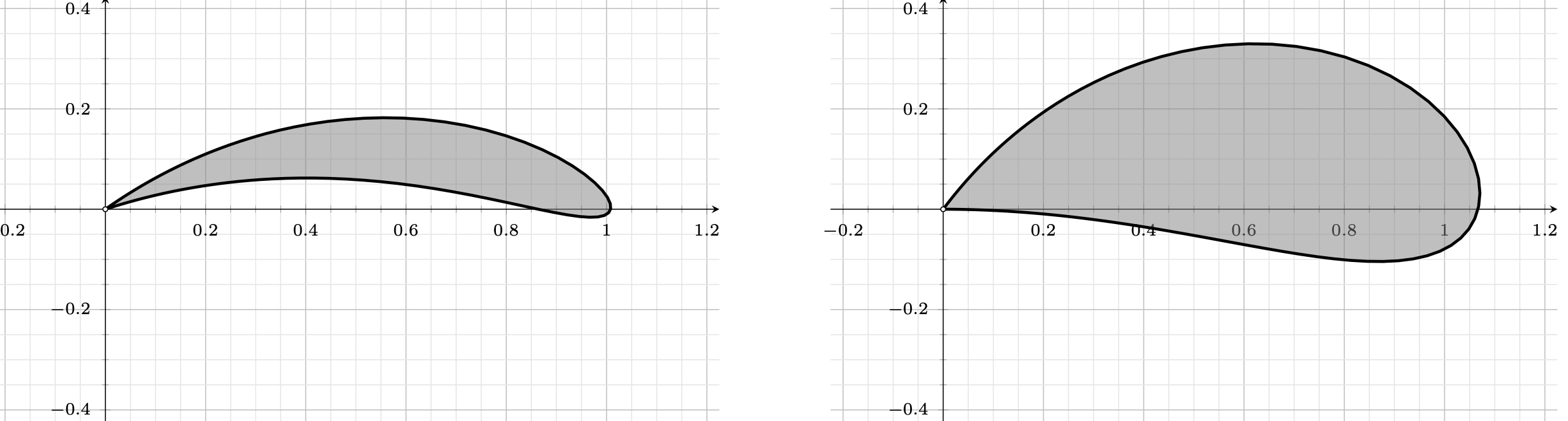}
    \caption{\small $B(C_\phi)$ for (i) $\phi(z) = \left(-\frac{1}{5} + \frac{3}{5}i\right)z + \frac{1}{5}$ (left) and (ii) $\phi(z) = \frac{i}{2}z + \frac{1}{2}$. Both ranges are not convex with the convexity condition~(\ref{ConvexityCondition}) failing: (i) $\frac{12\sqrt{11}}{25} \not\le  \frac{3}{25}$; (ii) $ 1 \not\le \frac{3}{8}$.}\label{ConFigs2}
\end{figure}

  Figures \ref{ConFigs1} and \ref{ConFigs2} below are computer-generated plots of $B(C_\phi)$ illustrating the application of the convexity condition (\ref{ConvexityCondition}). 
Note that for $\phi(z) = az + b$,  $B(C_\phi)$  depends only $|b|$: by (\ref{modbonly}), if $\theta_0$ is an argument of $b$, then 
$$
B(C_\phi) = \left\{  \frac{(1-r^2)}{1 - ar^2-|b|re^{-i(\theta-\theta_0)}}: 0\le r < 1,  \theta_0\le \theta < \theta_0+ 2\pi \right\}.
$$  For those Berezin ranges plotted in Figures \ref{ConFigs1} and \ref{ConFigs2}, note by Theorem~\ref{ADT} every boundary point of $B(C_\phi)$ other than $0$ belongs to $B(C_\phi)$ because $\phi$ is continuous on $\D^-$ and has no fixed point on $\partial \D$.  Also note that the Berezin range on the left in Figure~\ref{ConFigs1}, which is for the composition operator with symbol $\phi(z) = z/2 + 1/4$, $B(C_\phi)$ appears to be symmetric relative to the real axis.  

   \subsection{Symmetry of $B(C_\phi)$ relative to the real axis} Cowen and Felder ask \cite[Question 5.3]{cowen2022convexity}: For a composition operator $C_\phi$ acting on $H^2$, when is $B(C_\phi)$ symmetric about the real line?  It is easy to see that if $\phi$ is real on the real axis (which is true of $\phi(z) = z/2 + 1/4$), then $B(C_\phi)$ is symmetric relative to the real axis.  We provide a more general sufficient condition for symmetry below. 
   
  \begin{proposition}\label{symprop} Suppose that $\phi$ is an analytic selfmap of $\D$ such that $\phi$ maps a diameter of the open unit disk into itself; then $B(C_\phi)$ is symmetric relative to the real axis.
  \end{proposition}
  \begin{proof} Let $\phi$ map the diameter $\{te^{i\theta}: -1 < t < 1\}$ of the open unit disk into itself, where $\theta\in [0, 2\pi)$.  Then $\psi(z) := e^{-i\theta}\phi(e^{i\theta} z)$ is real on the real axis. Hence $\psi$ has real Maclaurin coefficients and $\overline{\psi(z)} = \psi(\bar{z})$ for every $z\in \D$.  Suppose that $w\in B(C_\psi)$ so that $w = \frac{1-|\alpha|^2}{1-\bar{\alpha}\psi(\alpha)}$ for some $\alpha \in \D$. Then
  $$
  \bar{w} =  \frac{1-|\alpha|^2}{1-{\alpha}\psi(\bar{\alpha})} = \widetilde{C_\psi}(\bar{\alpha}) \in B(C_\psi),
  $$
  and it follows that $B(C_\psi)$ is symmetric relative to the real axis.

  Let $\sigma(z) = e^{i\theta} z$. It is easy to check that $C_\sigma^*$ is the composition operator with symbol $z\mapsto e^{-i\theta}z$, so that $C_\psi = C_\sigma C_\phi C_\sigma^*$.  Also note that $C_\sigma^* k_\alpha = \sqrt{1-|\alpha|^2}K_{e^{i\theta}\alpha} = k_{e^{i\theta}\alpha}$, with the final equality holding because $|\alpha| = |e^{i\theta}\alpha|$. Thus,
   \begin{align*}
B(C_\psi)  & =  B(C_\sigma C_\phi C_\sigma^*) \\
& =  \{ \langle C_\phi C_\sigma^* k_{\alpha},  C_\sigma^*k_{\alpha} \rangle: \alpha \in \D\}\\
  &=  \{\langle C_\phi k_{e^{i\theta}\alpha}, k_{e^{i\theta}\alpha}\rangle: \alpha \in \D\}\\ 
  & = \{ \langle C_\phi k_{\alpha},  k_{\alpha} \rangle: \alpha \in \D\}\\
  & = B(C_\phi).
  \end{align*}
  Hence, $B(C_\phi)$ is also symmetric relative to the real axis. 
  \end{proof}

    \begin{remark}
    
    \begin{itemize}
    \item[(i)] The unitary composition operators on $H^2$ are precisely those whose symbols have the form $z\mapsto \zeta z$ for some $\zeta\in \partial \D$ (see, e.g., Theorem 8.2 and Exercise 8.1.3 of \cite{cowen1995composition}). Observe that the second paragraph of proof of the preceding proposition shows that   for any bounded linear operator $T$ on $H^2$, if $C_\sigma$ is unitary, then $B(C_\sigma T C_\sigma^*) = B(T)$.  More generally,  conjugation of $T$ by a unitary weighted composition operator leaves $B(T)$ unchanged \cite[Exercise 1, p.\ 160]{zhu2007operator}. 
    \item[(ii)]  Proposition 4.3 of \cite{cowen2022convexity} establishes that $B(C_\phi)$ is symmetric if $\phi$ is the automorphism of $\D$ given by $\phi(z) = (z -a)/(1-\bar{a}z)$, where $a\in \D$.  If $a\ne0$, then it is easy to check that this automorphism maps the diameter $\{t a/|a|: t\in (-1,1)\}$ to itself (and if $a = 0$, then $\phi(z) = z$ and $\phi$ maps every diameter to itself). Thus, the symmetry result of \cite[Proposition 4.2]{cowen2022convexity} also follows from Proposition~\ref{symprop} above.    
    \end{itemize}
    \end{remark}

    For selfmaps of $\D$ having the form $\phi(z) = az + b$, the converse of Proposition~\ref{symprop} holds.  
  Let $\phi(z)= az + b$ be a selfmap of $\D$ and suppose that $B(C_\phi)$ is symmetric relative to the real axis. To show that $\phi$ maps a diameter of $\D$ into itself it is sufficient to show $a$ is real, in which case  $\phi$ maps the interval $(-1,1)$ into itself if $b = 0$, and if $b\ne 0$ and $\theta$ is an argument of $b$,  then $\phi$ maps the diameter $\{te^{i\theta}: -1 < t < 1\}$ into itself.     
  
  We have already seen that $B(C_\phi)$ is symmetric relative to the real axis if $\phi$ fixes a point of $\partial \D$, and in this case $a$ is real ($0< a < 1$).  
  
  Suppose that $\phi$ does not fix a point of $\partial\D$ and $d = \sqrt{|1-a|^2 -|b|^2}$.
Because $B(C_\phi)$ is symmetric relative to the real axis, so is the set $1/B(C_\phi): = \left\{\frac{1}{z}: z\in B(C_\phi)\right\}$. We show in Section~\ref{Section4} that $1/B(C_\phi)$ has boundary equal to the hyperbola branch $\mathcal{h}: = \{w\in \CP: |w-a| - |w-1| = d\}$, which has  vertex $\frac{1+a}{2} + \frac{d}{2}\frac{1-a}{|1-a|}$ while the hyperbola of which $\mathcal{h}$ is a branch has foci $1$ and $a$.  The set  $1/B(C_\phi)$ consists of $\mathcal{h}$ and the component of $\CP\setminus\mathcal{h}$ containing the focus $1$.  (We remark that if  $b =0$,   $1/B(\phi)$ is the degenerate hyperbola consisting of the ray emanating from $1$ in the direction of the vector $1-a$.)  
Thus, the hyperbola-branch-bounded region $1/B(\phi)$ is symmetric relative to the real axis if and only if the focus $a$ is real, completing our proof of the converse of Proposition~\ref{symprop}  for selfmaps of $\D$ of the form $\phi(z) = az + b$.   Note that Figures \ref{ConFigs1} and \ref{ConFigs2} above illustrate this symmetry result: $B(C_\phi)$ is symmetric relative to the real axis if and only if the coefficient of $z$ for $\phi$ is real.

Another example of a composition operator whose Berezin range is symmetric relative to the real axis is the composition operator with symbol $\phi(z) = z/(2-iz)$, which maps the diameter $\{it: -1 < t < 1\}$ into itself.  Before providing a computer-generated image of   $B(C_\phi)$ for $\phi(z) = z/(2-iz)$, we note that  $\phi(z) = C_\sigma C_{\phi_2} C_\sigma^*$ where $\phi_2(z) := z/(2-z)$ and $\sigma$ is the selfmap of $\D$ given by $\sigma(z) = iz$; thus, $B(C_\phi) = B(C_{\phi_2})$ (see the proof of Proposition~\ref{symprop}). Figure~\ref{BerezinHorn} below is a plot of $B(C_\phi) = B(C_{\phi_2})$.  Now observe that the selfmap $\phi_2(z) = z/(2-z)$ is a member the family $\mathcal{F}_2$ (of (\ref{FamF2})). 

\subsection{$B(C_{\phi_s})$ for $\phi_s\in \mathcal{F}_2$}
We have promised to give a complete description of the Berezin ranges of all members of the family $\mathcal{F}_2$. Our method for obtaining the description depends upon  Lemma~\ref{TLL} and the idea of the first paragraph of the proof of Theorem~\ref{diskunion}.
 
 \begin{proposition}
\label{sobc} Suppose that $\psi$ is a nonconstant linear-fractional selfmap of $\D$, that $n$ is a positive integer, and that $\phi(z) = z\psi(z^n)$; then, for each $r \in [0, 1)$,   $\BTC(\{z: |z| = r\})$  is a circle {\rm(}degenerate if $r = 0${\rm)}, $B(C_\phi) = \cup_{r\in [0, 1)} \BTC(\{z: |z| = r\})$, and the interior of each of these image circles is contained in $B(C_\phi)$. 
 \end{proposition}
 \begin{proof}  The image of the circle $\{z: |z| = r\}$, $0 < r < 1$, under the Berezin transform of the composition operator $C_{\phi}$ is 
 $$
 f_r\left(\vstrut \{z: |z| = 1\}\right),
 $$ 
 where $f_r$ is the linear-fractional transformation $f_r(z) = \frac{1-r^2}{1-r^2\psi(r^nz)}$ (and the image circle is traversed $n$ times). That $B(C_\phi) = \cup_{r\in [0, 1)} \BTC(\{z: |z| = r\})$ is clear and that the interior of the circle $\BTC(\{z: |z| = r\})$ is contained in $B(C_\phi)$ follows from Lemma~\ref{TLL}.\end{proof}
 
Observe that $\phi_2(z) = z/(2-z) = z \psi(z)$, where $\psi$ is the linear-fractional selfmap $\psi(z) = 1/(2-z)$, so that the preceding corollary applies to $\phi_2$.   In Figure~\ref{BerezinHorn} below, we illustrate $B(C_{\phi_2})$ as a union of the circles $\widetilde{C_{\phi_2}}(\{z: |z| = r\})$, $0 \le r < 1$, as well as its symmetry relative to the real axis. 

   We now describe the method we use to completely characterize the Berezin range of $C_{\phi_2}$. 
For $0\le r < 1$, the image circles that determine the Berezin range of $C_{\phi_2}$ have diameters extending from 
$$
 \frac{1-r^2}{1-\frac{r^2}{2+r}} \quad  \text{to} \quad  \frac{1-r^2}{1-\frac{r^2}{2-r}},
 $$
 and thus, the centers are given by $c(r) = \frac{4-3r^2}{4-r^2}$, and radii by $\rho(r) = \frac{r^3}{4-r^2}$.  For $r\in [0,1)$, we let $(c(r), \rho(r))$ denote the circle with center $c(r)$ and radius $\rho(r)$.    It is easy to verify that as $r$ increases,  the centers move to the left and the radii increase.  Specifically, as $r$ increases from $0$ to $1$,  the centers $c(r)$ decrease from $1$ to $1/3$ and the radii $\rho(r)$ increase from $0$ to $1/3$.   Thus, in particular, the circles $(c(r),\rho(r))$, as $r\to1^-$, approach the {\it limit circle} having center $c(1) = 1/3$ and radius $\rho(1) = 1/3$ (which is the circle with dashed boundary in Figure~\ref{BerezinHorn}).   Observe that every number strictly within the limit circle belongs to $B(C_{\phi_2})$ by Proposition~\ref{sobc} above (each such point being inside $(c(r), \rho(r))$ for $r$ sufficiently close to $1$).   We determine the upper boundary of $B(C_{\phi_2})$ via a maximization process (and the lower boundary is simply the reflection in the real axis of the upper because of the symmetry of $B(C_{\phi_2})$ about the real axis).   Because $\phi_2(0) = 0$,  the Berezin radius of $C_{\phi_2}$ is $1$, which is attained at $1$. Thus, the only point in $B(C_{\phi_2})$ with real part 1 is $1$.
 
 The upper boundary of $B(C_{\phi_2})$ is the graph of the function $h: [0,1)\to \CP$ defined by 
 $$
 h(x):= \sup\{\Im(z): \Re(z) = x\ \&\ z\in B(C_{\phi_2})\},
 $$
 where the graph of $h$ is equivalent to $\{x+ih(x), 0\le x < 1\}$. We think of $h(x)$ as the height of the boundary curve directly above $x$.  
 Let $x\le 1/3$. We claim that $h(x)= \sqrt{1/9 - (x-1/3)^2}$, the height of the limit circle above $x$.  Because all points strictly inside the limit circle lie in $B(C_{\phi_2})$, we know $h(x)\ge  \sqrt{1/9 - (x-1/3)^2}$.  Also, no point in $\{x + i\kappa \sqrt{1/9 - (x-1/3)^2}: \kappa \ge 1\}$ can lie in $B(C_{\phi_2})$ because any point in $B(C_{\phi_2})$ would be on a circle with center $c(t_0) > 1/3$ and radius $\rho(t_0) < 1/3$ for some $t_0\in[0, 1)$, and if $z\in \{x + i\kappa\sqrt{1/9 - (x-1/3)^2}: \kappa \ge 1\}$, then $|z-c(t_0)| > |z-1/3| \ge 1/3$. Thus, $z$ cannot lie on the circle $(c(t_0),\rho(t_0))$. It follows that the upper boundary of $B(C_{\phi_2})$ over $[0, 1/3]$ is the graph of $h(x) =\sqrt{1/9 - (x-1/3)^2}$ and the graph is not in $B(C_{\phi_2})$.  
 
 Let $x\in (1/3, 1)$. Then $x = c(t_0)$ for some $t_0\in (0,1)$   We know that all points on and inside the circle  $(c(t_0), \rho(t_0))$ lie in $B(C_{\phi_2})$; thus $h(x) \ge \rho(t_0)$.  We claim that $h(x)$ exceeds $\rho(t_0)$.  A computation shows that for $\epsilon > 0$ sufficiently small with $t_0 + \epsilon < 1$, we have $\rho(t_0 + \epsilon)^2 - \left(c(t_0) - c(t_0+ \epsilon)\right)^2 > \rho(t_0)^2$; in other words, for such an $\epsilon$, the circle $(c(t_0 + \epsilon),\rho(t_0 + \epsilon))$ lies in $B(C_{\phi_2})$ and the point $P$ on this circle directly above $x$ (i.e., on the line $\{x+iy: y> 0\}$) lies above $x+i\rho(t_0)$. Note that because the interior of $(c(t_0 + \epsilon),\rho(t_0 + \epsilon))$ is in $B(C_{\phi_2})$, the entire segment $\{x + it\Im(P): t\in [0, 1]\}$ lies in $B(C_{\phi_2})$.
  It follows that 
\begin{equation}\label{bcphi2}
h(x) = \sup\{ \sqrt{  \rho(t)^2-  (x - c(t))^2}: t_0 \le t < 1\} \quad (\text{see Figure~\ref{MaxMethod})}, 
\end{equation}
and that $B(C_{\phi_2}) \cap \{z: \Re(z) = x\}$ is either $\{x + iy: -h(x) < y < h(x)\}$ or $\{x + iy: -h(x) \le y \le h(x)\}$.  The reason we can ignore $t< t_0$ in the preceding supremum of (\ref{bcphi2}) is that for any $t$ satisfying $0\le t < t_0$, $\rho(t) < \rho(t_0)$ which means $(c(t), \rho(t))$ either does not intersect $\{x + iy: y >0\}$ or intersects $\{x + iy: y >0\}$ below $x + i\rho(t_0)$.  
 
 For $x\in (1/3, 1)$, $c(t) = x$ has solution $t = 2\sqrt{(1-x)/(3-x)}$.  Thus, for a fixed $x\in (1/3, 1)$, we seek the supremum of the range of the function $f(t) =  \sqrt{  \rho(t)^2-  (x-c(t))^2}$ over the interval $ [2\sqrt{(1-x)/(3-x)}, 1)$.  To simplify the algebra for determining the range of $f$ and the corresponding supremum, we differentiate $$
 F(t) = f(t)^2 =   \rho(t)^2-  (x-c(t))^2 = t^2-x^2+6x-1 + \frac{16x}{t^2-4},
 $$ where  $c(t) = \frac{4-3t^2}{4-t^2}$ and  $\rho(t) = \frac{t^3}{4-t^2}$.  We have 
$$
 F'(t)  = t\left(\frac{2(4-t^2)^2 - 32x}{(4-t^2)^2}\right)  = t\left(\frac{\left(2\sqrt{1-\sqrt{x}} -t\right)\left(2\sqrt{1-\sqrt{x}} +t)\right)\left((4-t^2) + 4\sqrt{x}\right)}{(4-t^2)^2}\right).
 $$
 Because $x\in (1/3,1)$ we see that the sign $F'$ is determined by the sign of its factor $t\mapsto 2\sqrt{1-\sqrt{x}} -t$. Note that  $2\sqrt{1-\sqrt{x}} >  2\sqrt{(1-x)/(3-x)}$ for $x \in (1/3,1)$ and that $2\sqrt{1-\sqrt{x}} = 1$ when $x=9/16$ and is greater than $1$ when $x < 9/16$. 
 
  Let $x\in (1/3, 9/16]$.  We see that that $F'$ is positive on $[2\sqrt{(1-x)/(3-x)}, 1)$; thus, $F$ is increasing on this interval and $f$ is increasing on this interval too.  Hence,  the supremum of the range of $f(t)$  on $[2\sqrt{(1-x)/(3-x)}, 1)$ is $f(1)$ and the supremum is not attained. Thus for $x\in (1/3, 9/16]$, 
  $$
  h(x) = f(1) =  \sqrt{  \rho(1)^2-  (x-c(1))^2} = \sqrt{  1/9 -  (x- 1/3)^2}.
  $$
  We see that the upper boundary of $B(C_{\phi_2})$ coincides with upper boundary of the limit circle over $(1/3, 9/16]$ and that this portion of boundary is not contained in $B(C_{\phi_2})$.  

Finally, if $x >9/16$, then $2\sqrt{1-\sqrt{x}}$ is a critical number of $F$ in   $[2\sqrt{(1-x)/(3-x)}, 1)$ yielding the maximum value of $F$ on this interval. Thus, for $x\in (9/16, 1)$, we have 
$$
h(x) =  f\left( 2\sqrt{1-\sqrt{x}})\right) = (1-\sqrt{x})\sqrt{4 -(1+\sqrt{x})^2},
$$
so that $x + ih(x)$ lies on the circle $(c(2\sqrt{1-\sqrt{x}}\,), \rho(2\sqrt{1-\sqrt{x}}\,))$. Hence, the
graph of $h$ over $(9/16, 1]$  belongs to $B(C_{\phi_2})$ and is its upper boundary.  We have shown the upper boundary of $B(C_{\phi_2})$ is the graph of 
$$
h(x) = \left\{\begin{array}{cc} \sqrt{1/9 - (x-1/3)^2}  & \text{if}\ 0 \le x \le 9/16,\\ 
(1-\sqrt{x})\sqrt{4 -(1+\sqrt{x})^2} & \text{if}\ 9/16 < x \le 1;
\end{array} \right.
$$
 moreover, that portion of the graph of $h$ over $[0,9/16]$ is not contained in $B(C_{\phi_2})$ while that portion over $(9/16, 1]$ is contained in $B(C_{\phi_2})$.  

The process described above works for any member of the family $\mathcal{F}_2$. 
 Let $s>1$. The limit circle for $B(C_{\phi_s})$ has center $1/(1+s)$ and radius $1/(1+s)$. The limit circle provides the boundary of $B(C_{\phi_s})$ over the interval $[0, 1/s +1/(2s^2) -1/(2s^3)]$ and this portion of the boundary is not contained in $B(C_{\phi_s})$.  Finally, over the interval  $(1/s +1/(2s^2) -1/(2s^3), 1]$, the upper boundary of $B(C_{\phi_s})$ is the graph of 
$$
h(x) = \sqrt{-1 + 2 (s - 2 s^3 + s^4) + 2 x + 2 (-1 + s) s x - x^2 + (2 - 2 s)\sqrt{(-1 + s) s^3 (-2 - s + s^2 + 2 x)}},
$$
and the graph is contained in $B(C_{\phi_s})$.   Figure~\ref{Bphi5}  below displays $B(C_{\phi_5})$.

\begin{figure}[h]
    \centering
    \begin{minipage}{0.45\textwidth}
        \centering
        \includegraphics[width=0.92\textwidth]{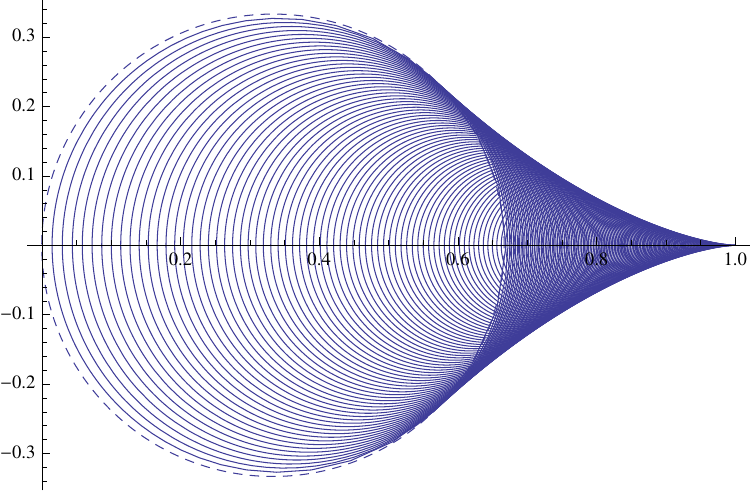} 
        \caption{\small $B(C_{\phi_2})$ as a union of circles.} \label{BerezinHorn}
    \end{minipage}\hfill
    \begin{minipage}{0.53\textwidth}
        \centering
        \includegraphics[width=\textwidth]{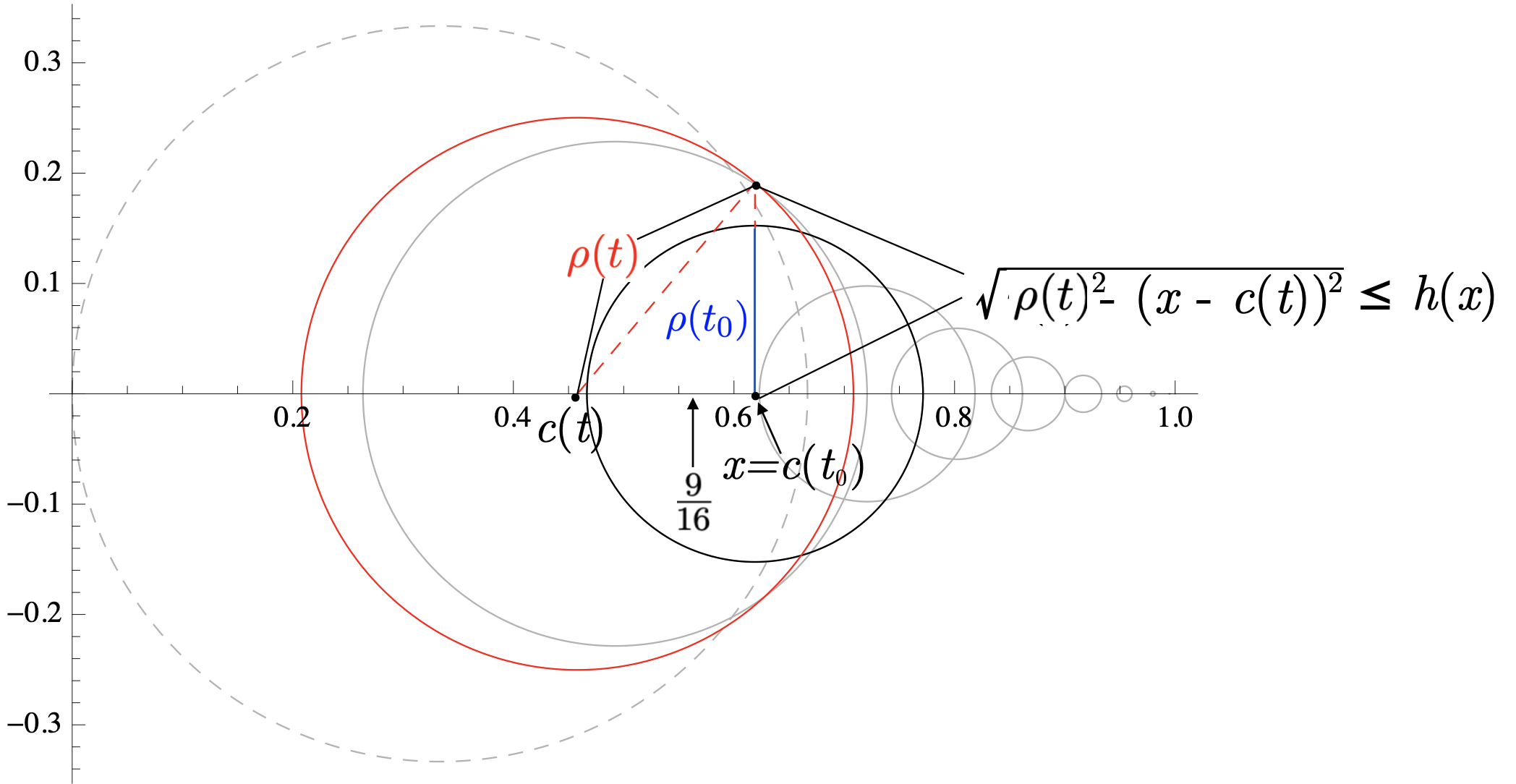} 
        \caption{\small $h(x) = \sup\{ \sqrt{  \rho(t)^2-  (c(t) - x)^2}: t_0 \le t < 1\}$.}\label{MaxMethod}
    \end{minipage}
    \end{figure}

    \begin{figure}[h]
    \centering
     \includegraphics[height=1.9in]{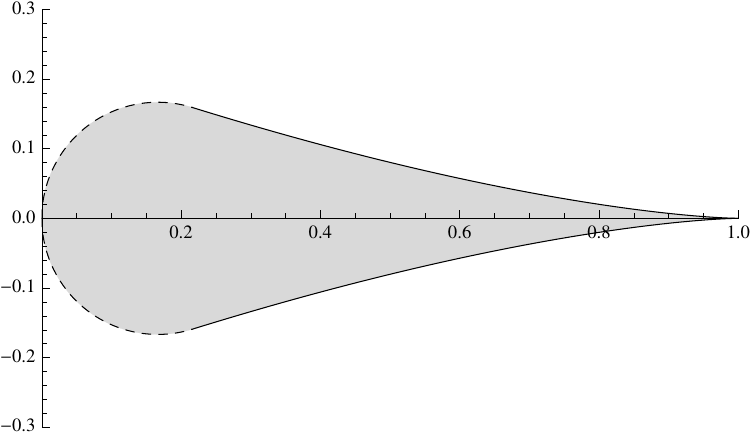} 
        \caption{\small $B(C_{\phi_s})$ for $s = 5$.}\label{Bphi5}
    \end{figure}

\subsection{Corner Points}

Observe that the origin appears to be a ``corner point'' of all the Berezin ranges of Figures~\ref{ConFigs1}  and \ref{ConFigs2}, which are for selfmaps of the form $\phi(z) = az + b$ that do not fix a point on $\partial \D$.   In contrast, $1$ appears to be a corner point for the Berezin ranges of Figures~\ref{BerezinHorn} and \ref{Bphi5}, which are for the selfmaps $\phi_2$ and $\phi_5$ of the family $\mathcal{F}_2$.  

We use triangles to define ``corner point.''  We assume that a triangle $T$ consists of its three-sided boundary as well as its interior.  For a point $p$ in $\CP$ and a positive number $r$, we let $D(p; r)$ be the open disk centered at $p$ having radius $r$.

\begin{definition}\label{cornerdef} A point $p$ on the boundary of $B(C_\phi)$ is a corner point of $B(C_\phi)$ provided there is a triangle $T$ with vertex $p$ and a number $\epsilon > 0$ such that  $B(C_\phi) \cap D(p; \epsilon) \subseteq T$.
\end{definition}

\begin{proposition} Suppose that $\phi$ is a selfmap of $\D$ satisfying the hypotheses of either Theorem~\ref{diskunion} or Proposition~\ref{sobc}. If  $p$ is a corner point of $B(C_\phi)$, then $p = 0$ or $p =1$.
\end{proposition}
\begin{proof} Let $p$ be a corner point of $B(\C)$.  Let $T$ be a triangle with vertex $p$ and let $\epsilon > 0$ be such that $B(C_\phi) \cap D(p; \epsilon) \subseteq T$.  Let $\theta$ be the angle of $T$ at vertex $p$ (necessarily $\theta < \pi$), and observe that if $C$ is a circle in $T$ having radius $s >0$,  then $\min\{|z-p|: z\in C\} \ge s(1/\sin(\theta/2)  -1 )$. 

Let $f$ be the Berezin transform of $C_\phi$. For $0 \le r < 1$,  let $\beta_r$ be the circle $\{f(re^{it}): 0 \le t < 2\pi\}$ and let $F_r$ be the closed disk consisting of $\beta_r$ and its interior.   Note $F_r\subset B(C_\phi)$ because $\phi$ satisfies the hypotheses of either Theorem~\ref{diskunion} or Proposition~\ref{sobc}.    Let $\rho(r)$ be the radius of $\beta_r$ so that $\rho$ is continuous on $[0, 1)$ and positive on $(0, 1)$. 

 Because $p$ is a boundary point of $B(C_\phi)$, there is a sequence $(z_n)$ of points in $B(C_\phi)$ such that $(z_n)$ converges to $p$. For each $n\in \N$, chose $\beta_{r_n}$ such that $\beta_{r_n}$ contains the point $z_n$.   Claim: $(\rho(r_n))$, the sequence of radii of the circles $(\beta_{r_n})$, converges to $0$. Suppose, in order to obtain a contradiction, that there is a subsequence $(r_{j_n})$ of $(r_n)$ and a number $\gamma > 0$ such that $\rho(r_{j_n}) >\gamma$ for all $n\in \N$.  For each $n$, let $\tau_n$ be a circle contained in $F_{r_{j_n}}$, with diameter equal to $\min\{\epsilon/2, \gamma\}$  that is internally tangent to $\beta_{r_{j_n}}$ at $z_{j_n }$.  Because $(z_{j_n})$ converges to $p$, there is an $N\in \N$ such that for all $n\ge N$, $|z_{j_n} - p| < \epsilon/2$.  Let $n\ge N$ be arbitrary. Because  the diameter of $\tau_{n}$ is at most $\epsilon/2$ and $\tau_{n}$ contains $z_{j_n}$, whose distance to $p$ is less than $\epsilon/2$, every point on $\tau_{n}$ lies within  $\epsilon$ distance of $p$; i.e., $\tau_{n} \subset D(p;\epsilon)$. Also,  $\tau_{n}\subset  F_{r_{j_n}}\subset B(C_\phi)$.  Thus, $\tau_{n} \subset B(C_\phi)\cap D(p;\epsilon)\subseteq T$.  We conclude that 
 $$
 \min\{|z-p|: z\in \tau_{n}\} \ge \min\{\epsilon/4, \gamma/2\}(1/\sin(\theta/2) -1) >0
 $$
 which, because $n\ge N$ is arbitrary, implies $|z_{j_n}- p| \ge \min\{\epsilon/4, \gamma/2\}(1/\sin(\theta/2) -1)$ for all $n\ge N$, a contradiction.  Thus, $(\rho(r_n))$, the sequence of radii of the circles $(\beta_{r_n})$, converges to $0$.

 Because $r_n \in [0, 1)$ for each $n$, $(r_n)$ has a convergent subsequence $(r_{j_n})$. Note that $(r_{j_n})$ must converge to either $0$ or $1$; otherwise, it converges to a number $r$ in $(0,1)$ so that $(\rho(r_n))$ converges to $\rho(r) > 0$, a contradiction.

 Suppose  $(r_{j_n})$ has limit $1$.  By the proof of Theorem 3.1 of \cite{bourdon2002zero}, there is a $t_0\in [0, 2\pi)$ such that 
 $f(r_{j_n}e^{it_0}) \to 0$.  Because,  for each $n\in \N$, both  $f(r_{j_n}e^{it_0})$ and $z_{j_n}$ lie on the circle  $\beta_{r_n}$, and the sequence of corresponding radii $(\rho(r_{j_n}))$ converges to $0$, we see that $(|z_{j_n} - f(r_{j_n})e^{it_0}|)$  converges to $0$. Hence, $(z_{j_n})$ converges to $0$, so that $p = 0$.
 
 Suppose  $(r_{j_n})$ has limit $0$. Because for each $n\in \N$, $\beta_{r_{j_n}} =  \{f(r_{j_n}e^{it}): 0 \le t < 2\pi\}$, every sequence whose $n$-th term lies in $\beta_{r_{j_n}}$, e.g., $(z_{j_n})$ must have limit $f(0) = 1$ by the continuity of $f$ (on $\D$). Thus, $p =1$.  
 \end{proof}
  
\begin{figure}[h]
    \centering
     \includegraphics[height=1.3in]{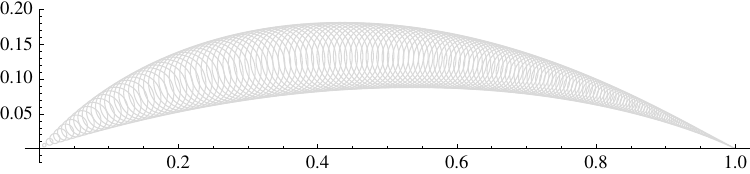} 
        \caption{\small$B(C_{\phi})$ for $\phi(z) = iz/(2-z)$ has 0 and 1 as corner points by Proposition~\ref{0CPTProp} and Theorem~\ref{1CPThm}.} \label{2CPexample}
    \end{figure}

   As an application of Proposition~\ref{0CPTProp} below, we show that $0$ is a corner point of $B(C_\phi)$ for the selfmap $\phi(z) = az + b$ if $b\ne 0$ and $\phi$ does not fix a point of $\partial \D$. Parts (i) and (ii) of Theorem~\ref{T:Berezin_range_az_b} below show that if $\phi(z) = az$ is a selfmap of $\D$, then, again, $0$ is a corner point of $B(\C)$.   We will also show that $1$ is a corner point of $B(C_\phi)$ if and only if $\phi(0) = 0$ (a property that all selfmaps $\phi_s\in \mathcal{F}_2$ share). Figure~\ref{2CPexample} below provides an interesting example of a Berezin range having both $0$ and $1$ as corner points.

 \begin{proposition}\label{0CPTProp}  Let $\phi$ be an analytic selfmap of $\D$.  Suppose that there is a positive constant $\gamma$ such that 
\begin{equation}\label{0CPT}
  \Re(1-\bar{\alpha}\phi(\alpha)) \ge\gamma \ \text{for all}\ \alpha \in \D;
\end{equation}
  then $0$ is a corner point of $B(C_\phi)$.
  \end{proposition}
  \begin{proof} 
  Let $z$ be in $B(C_\phi)$, so that there is some $\alpha \in \D$ for which $\widetilde{C_\phi}(\alpha) = z$. Because $z$ lies in the right halfplane $z = |z|e^{i\theta}$ for some $\theta \in (-\pi/2, \pi/2)$. We show that   $\theta \in (-\tan^{-1}(1/\gamma), \tan^{-1}(1/\gamma)$.  We have
  $$
 \frac{1}{z} = \frac{1}{1-|\alpha|^2} \left(1-\bar{\alpha}\phi(\alpha)\right).  
 $$
 Thus, $-\theta$ is an argument of  $\left(1-\bar{\alpha}\phi(\alpha)\right)$, a number whose real part is $\ge \gamma >0$ and  whose imaginary part  lies in $(-1, 1)$ because  $\phi$ is a selfmap of $\D$.   Thus, $\tan(\theta)$ lies in $(-1/\gamma, 1/\gamma)$, as claimed.  We see that $z$ lies in the sector determined by the rays emanating from the origin making angles $-\tan^{-1}(1/\gamma)$ and $\tan^{-1}(1/\gamma)$ with the positive real axis, respectively.  Because $z\in B(C_\phi)$ is arbitrary, $B(C_\phi)$ lies within this sector and $0$ is thus a corner point of $B(C_\phi)$.
 \end{proof}  

 Applications of Proposition~\ref{0CPTProp}:
 \begin{itemize}
     \item[(i)]
  Let $\phi(z) = iz/(2-z)$. For each $r \in [0,1)$,   $\{1-\bar{\alpha}\phi(\alpha): |\alpha| = r\} = \{1-if_r(\{z: |z| = 1\})\}$, where $f_r(z) = r^2/(2-rz)$. Note that for $r\in (0,1)$, $f_r(\{z: |z| = 1\})$  is a circle in the open right halfplane with center on the real axis and radius $r^3/(4-r^2)$, while $f_0(\{z: |z| = 1\}) = \{0\}$, a the degenerate circle.  Thus, for each $r \in [0,1)$, $\{1-\bar{\alpha}\phi(\alpha): |\alpha| = r\}$ is a circle with center on the ray  $\{1-it: t\ge 0\}$ and with radius $r^3/(4-r^2)$. Because the radii increase to 1/3 as $r\to 1^-$, we see that condition (\ref{0CPT}) holds with $\gamma = 2/3$. Thus, $B(C_\phi)$ has $0$ as a corner point.  
  \item[(ii)]  Any $\phi$ such that $\|\phi\|_\infty < 1$ has $0$ as a corner point of $B(C_\phi)$. 
  \item[(iii)]  Suppose that $\phi(z) = az + b$ is a selfmap of $\D$ such that $b\ne 0$ and $\phi$ doesn't fix a point in $\partial \D$. Because, $\phi$ is a selfmap of $\D$, we have $|a| + |b|\le 1$ and because $b\ne 0$, we have $|a| < 1$.     Then, assuming $\Re(a) \ge 0$, we have
 $$
 \Re\left(1-\bar{\alpha}\phi(\alpha)\right) = \Re(1-|\alpha|^2a - b\bar{\alpha}) \ge 1 - \Re(a) -|b| \ge 1 -|a| -|b| \ge 0.
 $$
 We claim that $1 - \Re(a) -|b| > 0$; otherwise, the final two inequalities are equalities, and we have both $\Re(a) =|a|$ and   $|a| + |b| = 1$. In this case we have $a= |a|$ and $a+ |b| =1$.  However, $a+ |b| =1$ puts the fixed point $b/(1-a)$ of $\phi$ on $\partial \D$, contradicting our hypothesis.   Thus, if $\Re(a) \ge 0$,  then $\phi$  satisfies (\ref{0CPT}) with $\gamma =  1 - \Re(a) -|b|$.   If $\Re(a) < 0$, then $\phi$ satisfies (\ref{0CPT}) as well with $\gamma = 1 -|b|$.  Thus, by Proposition~\ref{0CPTProp}, $B(C_\phi)$ has $0$ as a corner point. 
 \end{itemize}  

 \begin{lemma}\label{BSL}  Suppose that $\phi$ is an analytic selfmap of $\D$ such that $\phi(0) = 0$, and $\phi$ is not an automorphism of $\D$ $($i.e., not a rotation of $\D)$; then, for every $\epsilon$ satisfying $0 < \epsilon < 1$, 
 \begin{equation}
   \sup \{|\BTC(\alpha)| : \epsilon \le  |\alpha| < 1\}  < 1.
\end{equation}
 \end{lemma}
 \begin{proof} Applying the generalized Schwarz Lemma stated immediately before Proposition~\ref{brgub} to the selfmap $z\mapsto \phi(z)/z$ if $z\ne 0$ and  $0\mapsto \phi'(0)$, we have, for all $z\in \D$,
 $$
 |\phi(z)| \le |z| \frac{|z| + |\phi'(0)|}{1 + |z||\phi'(0)|}.
 $$
 
 Let $\epsilon\in (0, 1)$.  We have, for each $\alpha \in \D$ satisfying $|\alpha| \ge \epsilon$,
  \begin{align*}
 |\BTC(\alpha)| &\le \frac{1-|\alpha|^2}{1-|\alpha|^2\frac{|\alpha| + |\phi'(0)|}{1 + |\alpha||\phi'(0)|}}\\
  & = 1 - \left( 1 - \frac{(1 + |\alpha|)(1 + |\alpha||\phi'(0)|)}{1+ |\alpha| + |\alpha|^2 + |\phi'(0)||\alpha|}\right)\\
   & = 1 - \frac{|\alpha|^2(1- |\phi'(0)|)}{1+ |\alpha| + |\alpha|^2 + |\phi'(0)||\alpha|}\\
   & \le 1 - \frac{|\epsilon|^2(1-|\phi'(0)|)}{3 + |\phi'(0)|}\\
   & < 1,
     \end{align*}
     where the final inequality holds because $|\phi'(0)|< 1$ by Schwarz's Lemma (because $\phi$ is not an automorphism of $\D$).
 \end{proof}            

  \begin{theorem}\label{1CPThm}  Suppose that $\phi$ is an analytic selfmap of $\D$; then $1$ is a corner point of $B(C_\phi)$ if and only if $\phi(0) = 0$.
 \end{theorem}
 \begin{proof}  
Let $\phi(0) = 0$.  If $\phi$ is an automorphism of $\D$, then $\phi(z) = \zeta z$ for some unimodular constant $\zeta$; and $B(C_\phi)$ is either the singleton set $\{1\}$ if $\zeta = 1$, or, if $\zeta \ne 1$, is, by parts (i) and (ii) of Theorem~\ref{T:Berezin_range_az_b} below, a half-closed circular arc $\mathcal{A}$ (possibly degenerate) passing through $1, 1/(2-\zeta)$, and $0$, where $\mathcal{A}$ contains the endpoint $1$ but not the endpoint $0$. Clearly, $\{1\}$ or $\mathcal{A}$ is contained in a triangle having $1$ as a vertex.
 
 Suppose that $\phi$ is not an automorphism (and $\phi(0) = 0$).  Throughout the remainder of the proof,   {\it argument} means ``principal branch of the argument,'' which we assume has range $[-\pi, \pi)$. Thus, the argument of a nonzero complex number is uniquely determined and we denote it by $\Arg$.   
 
 For $\alpha \in \D$, we have
 $$
 1- \BTC(\alpha)  = 1 - \frac{1-|\alpha|^2}{1-\bar{\alpha}\phi(\alpha)} 
   = \frac{|\alpha|^2\left(1- \frac{\phi(\alpha)}{\alpha}\right)}{1-\bar{\alpha}\phi(\alpha)}
  $$
  Thus, as $|\alpha| \to 0$, the argument of  $1-\BTC(\alpha)$ approaches the argument of $1-\phi'(0)$.  Let $\theta_0 = \Arg(1-\phi'(0))$.  Because $\phi$ is not an automorphism, we have, by Schwarz' s Lemma, $\phi'(0)\in \D$, so that $\theta_0\in (-\pi/2, \pi/2)$.  
  
  Let $\delta > 0$ satisfy $\delta < \min\{\pi/2 - \theta_0, \theta_0  + \pi/2\}$. 
    Choose $\epsilon \in (0, 1)$ small enough so that if $|\alpha| < \epsilon$, then $\Arg(1-\BTC(\alpha)) \in (\theta_0 - \delta, \theta_0 + \delta)$, which is an interval within $(-\pi/2, \pi/2)$. Let $\mathcal{R}_1$ be the ray emanating from $1$ making angle $-\theta_0+\delta$  with respect to the negative real direction and let $\mathcal{R}_2$ be the ray emanating from $1$ making angle $-\theta_0-\delta$ with respect to the negative real direction. Let $P_1$ and $P_2$ be the points where the rays $\mathcal{R}_1$ and $\mathcal{R}_2$, respectively, intersect the imaginary axis and let $T$ be the triangle with vertices $P_1$, $P_2$ and $1$.  
     
      Observe that we have constructed $T$ so that if $|\alpha| < \epsilon$, then $\BTC(\alpha)$ lies in $T$.  Let 
      $\beta = \inf\{|1-\BTC(\alpha)|: \epsilon \le |\alpha| < 1\}$ and observe $\beta \ge  1 -  \sup \{|\BTC(\alpha)| : \epsilon \le  |\alpha| < 1\} $. Thus, $\beta > 0$ by Lemma~\ref{BSL}.  Observe that $D(1; \beta)\cap B(C_\phi)\subseteq T$, which completes our proof that $1$ is a corner point of $B(C_\phi)$.  
      
      We complete the proof of the theorem by showing that if $\phi(0)\ne 0$, then $1$ is not a corner point of $B(C_\phi)$. Suppose that $\phi(0) \ne 0$ and let $\zeta = e^{i\theta_0}$, where $\theta_0 = \Arg(\phi(0))$.  Define the curve $s: (-1,1)\to \CP$ by $s(t) = \frac{1-t^2}{\rule{0in}{0.12in} 1-t\bar{\zeta}\phi(t\zeta)}$, so that $s(0) = 1$ and its image is contained in $B(C_\phi)$.  Now observe that $s'(0) = |\phi(0)| > 0$.  Thus, the image of $s$ has a horizontal tangent at $1$, and thus $1$ cannot be a corner point of $B(C_\phi)$.  
 \end{proof}
 
\begin{remark}
 \begin{itemize}
 \item[(i)] In the proof of the preceding theorem, the positive number $\delta$  that determines the angle of the triangle $T$ at its vertex $1$ may be taken to be as small as desired; thus,  when $\phi(0) = 0$, the corner point of $B(C_\phi)$ at $1$ is always a ``cusp point.''  
 \item[(ii)]  The number $0$ is seldom a corner point for the numerical range $W(C_\phi)$ of a composition operator $C_\phi$ because $0$ is typically in the interior of $W(C_\phi)$ (\cite[Theorem 5.3 \& Section 8]{bourdon2002zero}. In contrast, Proposition~\ref{0CPTProp} above and the applications discussed after it show that it is not uncommon for $0$ to be a corner point of $B(C_\phi)$.  However, if $\phi(0) = 0$, then $W(C_\phi)$ typically has $1$ as a corner point: in fact, by \cite[Theorem~7.3]{bourdon2002zero}, $1$ is a corner point of $W(C_\phi)$ if $\phi(0) = 0$ and $\phi$ is not an inner function.  Thus, for non-inner selfmaps $\phi$ of $\D$ fixing $0$, we see that $1$ is a corner point of $B(C_\phi)$ because it is a corner point of $W(C_\phi)$.  However, $1$ is not a ``cusp point'' of $W(C_\phi)$ when $\phi(0)= 0$ unless $\phi(z) = tz$ for some $t\in [-1,1]$; see  Theorem 3.8 of \cite{bourdon2002zero} and the discussion of when $1$ is a corner point of $W(C_\phi)$ on page 437 of \cite{bourdon2002zero}.
 \end{itemize}   
  \end{remark}

\section{A necessary and sufficient condition for convexity of $B(\C)$ for $\ph(z)=az+b$}\label{Section4}

We begin with a description of $B(C_\phi)$ for $\phi(z) = az +b$. While we have already established that $B(\C)$ is an open disk when $b\ne 0$ and $\phi$ fixes a point in $\partial \D$ (Corollary~\ref{ab1}), here we take a different  approach that also yields a description of $B(C_\phi)$ when $b\ne 0$ and $\phi$ does not fix a point on $\partial D$.  In addition, we describe $B(C_\phi)$ when $b =0$.

\begin{theorem}\label{T:Berezin_range_az_b}
Consider the operator $C_\phi:H^2 \rightarrow H^2$, where $\phi(z) = az+b$ is a selfmap of the disk, which is not the identity map.
\begin{itemize}
\item[(i)] If $b = 0$ and $a\in [-1,1)$, then $B(C_{\phi})$ is the half-closed line segment $(0,1]$.

\item[(ii)] If $b = 0$, $|a|\leq 1$, and $\Im(a)\ne 0$, then $B(C_{\phi})$ is the half-closed circular arc through the points $1,1/(2-a)$ and $0$, containing the endpoint $1$ but not the endpoint $0$.

\item[(iii)] If $b\neq 0$ and the fixed point of $\phi$ lies in $\partial\D$, then $0<a<1$ and $B(C_{\phi})$ equals the open disk centered at $1/(1+a)$ with radius $1/(1+a)$.

\item[(iv)] If $b\neq 0$ and the fixed point of $\phi$ lies in $\D$, then $B(C_{\phi})$ is the set of all $z\in\mathbb{C}\backslash\{0\}$ that satisfy the inequality
\begin{equation}
\label{Eqn:Berezin_range_az_b}
 \Big|\frac{1}{z}-a\Big| - \Big|\frac{1}{z}-1\Big| \ge d,
\end{equation}
where $d = \sqrt{|1-a|^2 - |b|^2}$.
\end{itemize}
\end{theorem}

\begin{proof}
As in the proof of Theorem \ref{diskunion}, the Berezin range is given by
\begin{align}
\label{Eqn:defn_Ber_range}
B(C_{\phi}) = \Big\{\frac{1-r^2}{1-ar^2-bre^{-i\theta}}: 0\leq r < 1, 0\leq \theta < 2\pi\Big\}.
\end{align}
Consider first the case $b = 0$. Then
\begin{align*}
B(C_{\phi})
& = \Big\{\frac{1-r^2}{1-ar^2}: 0\leq r < 1, 0\leq \theta < 2\pi\Big\}\\
& = \Big\{\frac{1-t}{1-at}: 0\leq t < 1\Big\}\\
& = \Phi([0,1)),
\end{align*}
where $\Phi(w) = \frac{1-w}{1-aw}$ is a linear fractional transformation. Note that $\Phi(0) = 1$, $\Phi(\frac{1}{2}) = \frac{1}{2-a}$, and $\Phi(1) = 0$ since $a\neq 1$. As a result, $\Phi([0,1))$ is the half-closed (possibly degenerated) circular arc through the points $1, \frac{1}{2-a}$ and $0$, containing the endpoint $1$ but not the other endpoint $0$. In the case $a$ is real, the arc degenerates into the half-closed interval $(0,1]$. We have thus proved $(i)$ and $(ii)$.

For the remainder of the proof, assume that $b\neq 0$. Equation \eqref{Eqn:defn_Ber_range} shows that a complex number $z$ belongs to $B(C_{\phi})$ if and only if there exists $r\in [0,1)$ and $\theta\in [0,2\pi)$ such that
\[\frac{1-r^2}{1-ar^2 - bre^{-i\theta}} = z,\]
which is equivalent to 
\begin{equation}
\label{Eqn:z_in_range}
(1-r^2) - (1-ar^2)z  = - bre^{-i\theta}z.
\end{equation}
The existence of such $r$ and $\theta$ solving \eqref{Eqn:z_in_range} is equivalent to the existence of a solution $r\in [0,1)$ to the equation
\[\Big|(1-r^2) - (1-ar^2)z\Big|^2 = |b|^2r^2|z|^2.\]
Setting $u = 1 - r^2$, we have $u\in (0,1]$ and
\[\Big|(1-az)u - (1-a)z\Big|^2 = |b|^2|z|^2(1-u).\]
Expanding the left side and grouping appropriate terms, we obtain the equation
\begin{equation}
\label{Eqn:u}
  |1-az|^2\cdot u^2 - \big(2\Re((1-a)z(1-\bar{a}\bar{z}))-|b|^2|z|^2\big)u +(|1-a|^2-|b|^2)|z|^2 = 0.  
\end{equation}
Applying the formula
\[2\Re(A\overline{B}) = |A|^2 + |B|^2 - |A-B|^2,\] with $A = (1-a)z$ and $B = 1-az$, we write
\begin{align*}
2\Re((1-a)z(1-\bar{a}\bar{z}))-|b|^2|z|^2 
& = 
   |(1-a)z|^2 + |1-az|^2 - |z-1|^2 - |b|^2|z|^2\\
& = d^2|z|^2 + |1-az|^2 - |z-1|^2,
\end{align*}
where $d = \sqrt{|1-a|^2-|b|^2}$. Equation \eqref{Eqn:u} then becomes
\begin{equation} 
\label{Eqn:u_simplified}
|1-az|^2\cdot u^2 - \big(d^2|z|^2 + |1-az|^2 - |z-1|^2\big)u + d^2|z|^2 = 0.
\end{equation}
We now divide the remaining argument into two cases.
\medskip

\underline{Case 1:} $d = 0$. In this case, $|b| = |1-a|$, so that the fixed point $b/(1-a)$ of $\phi(z) = az + b$ lies in $\partial \D$. By Proposition~\ref{LSMBFP},  $0<a\leq 1$; in fact, $0<a<1$ because we are assuming that $\phi$ is not the identity map. Equation \eqref{Eqn:u_simplified} now becomes
\[|1-az|^2\cdot u^2 - (|1-az|^2 - |z-1|^2)u = 0.\]
This equation has a root $0<u\leq 1$ if and only if $|1-az|^2 - |z-1|^2 > 0$, or equivalently,
$$(1-a^2)|z|^2 - 2(1-a)\Re(z) < 0\ \Longleftrightarrow (1+a)|z|^2 -2\Re(z)< 0\ \Longleftrightarrow  \Big|z - \frac{1}{1+a}\Big|^2 < \frac{1}{(1+a)^2},$$
where $1-a>0$ yields the first equivalence of the preceding line and where the final inequality means precisely that $z$ belongs to the open disk centered at $1/(1+a)$ with radius $1/(1+a)$. This proves $(iii)$.

\underline{Case 2:} $d>0$. In this case, the fixed point of $\phi$ lies in $\D$ and the necessary conditions for \eqref{Eqn:u_simplified} to have a root in $(0,1]$ are that $z\neq 0$ and
\begin{equation}
\label{Eqn:z_necessary_condition}
d^2|z|^2 + |1-az|^2 - |z-1|^2 \geq 2d|z|\cdot|1-az|,
\end{equation}
which is equivalent to
\begin{align*}
    0 & \leq  d^2|z|^2 + |1-az|^2 - |z-1|^2 -2d|z|\cdot|1-az|\\
    & = \big(d|z| - |1-az|\big)^2 - |z-1|^2\\
    & = \big(d|z| - |1-az| - |z-1|\big)\big(d|z| - |1-az| + |z-1|\big).
\end{align*}
Using the triangle inequality and the fact that $d <|1-a|$ (since $b\neq 0$) and $z\neq 0$, we estimate the first factor as follows
\begin{align*}
    d|z| - |1-az| - |z-1| & \leq d|z| -  |(z-1)+(1-az)| \\
    &= (d-|1-a|)|z| < 0.
\end{align*}
We then conclude that \eqref{Eqn:z_necessary_condition} is equivalent to
\begin{equation}
\label{Eqn:z_necessary_condition_simplified}
d|z| - |1-az| + |z-1| \leq 0,
\end{equation}
which, in turn, is equivalent to \eqref{Eqn:Berezin_range_az_b} since $z\neq 0$.

To complete the proof, we will show that Condition~(\ref{Eqn:z_necessary_condition_simplified}) is also sufficient for the existence of a root $u\in (0,1]$ to the equation \eqref{Eqn:u_simplified}. Let $z\neq 0$ satisfy \eqref{Eqn:z_necessary_condition_simplified}. In particular, $d|z|\leq |1-az|$. Because $d > 0$ and $|z|>0$, we see that $|1-az| > 0$, which together with \eqref{Eqn:z_necessary_condition}, shows that the quadratic equation \eqref{Eqn:u_simplified} has two positive roots. Since
\[d^2|z|^2 + |1-az|^2 - |z-1|^2\leq 2|1-az|^2\]
because $d|z| \leq |1-az|$, we conclude that at least one of the roots must belong to the interval $(0,1]$.
This proves $(iv)$ and completes the proof of the theorem.
\end{proof}

We note that parts (i) and (ii) of the preceding theorem relate to Theorem 4.1 of 
\cite{cowen2022convexity}, which states that  $B(C_\phi)$ is convex for a rotation-selfmap $\phi(z) = \zeta z$, $\zeta\in \partial \D$, if and only if $\zeta =1$ or $\zeta = -1$, with the proof indicating that when $\zeta = -1$, $B(C_\phi) = (0, 1]$.

We now obtain the following  necessary and sufficient condition for convexity of $B(C_\phi)$ when
$\phi(z) = az + b$:   $B(C_\phi)$ is convex  if and only if   $2d|\Im(a)| \le |b|(1-|a|^2)$, where  $d = \sqrt{|1-a|^2 - |b|^2}$.
It is easy to check the
validity of the preceding condition if $b = 0$ (encompassing the possibilities of cases (i) and (ii) of Theorem 4.1 as well as $\phi(z) = z$). Moreover, establishing its validity relative to case (iii) of Theorem~\ref{T:Berezin_range_az_b} is straightforward because case (iii)'s hypothesis on $\phi$ implies $C_\phi$ is convex.  Thus, to complete our convexity characterization, we focus on case (iv) of Theorem~\ref{T:Berezin_range_az_b}, where $B(C_\phi)$ is characterized as the set of reciprocals of points in the set 
   $$
\left\{z\in \CP: \big|z-a\big| - \big|z-1\big| \ge d\right\},
$$
whose boundary is a branch of a hyperbola in $\CP$ with foci $a$ and $1$, specifically, the branch having vertex $a + \frac{d}{2}(1-a)/|1-a|$.  Thus, to obtain our convexity characterization for $B(C_\phi)$, we begin with a lemma of independent interest characterizing when the convexity of a hyperbola branch (of a hyperbola in ``standard position'') is preserved by an inversion mapping. 

\begin{lemma}
\label{L:hyperbola_inversion}
Let $0<\delta<1$ and define
\begin{equation}\label{HBX}
\mathcal{X}_{\delta} = \Big\{w\in\mathbb{C}: |w+1| - |w-1| \geq 2\delta\Big\}.
\end{equation}
Let $w_0\in\mathbb{C}\setminus \mathcal{X}_\delta$. Then the image of $\mathcal{X}_\delta$ under the map $w\mapsto\frac{1}{w-w_0}$ is a convex set if and only if \[|\Im(w_0)|+\frac{\sqrt{1-\delta^2}}{\delta}\Re(w_0)\leq 0.\]
\end{lemma}

\begin{proof}
We see that the boundary of $\mathcal{X}_\delta$ is the right branch of the hyperbola with foci $-1$ and $1$, and constant difference $2\delta$. Some algebra (starting with squaring both sides of the inequality defining $\mathcal{X}_\delta$) yields
\[\mathcal{X}_\delta = \Big\{w\in\mathbb{C}: \frac{(\Re(w))^2}{\delta^2} - \frac{(\Im(w))^2}{1-\delta^2} \geq 1,\quad \Re(w) > 0\Big\}.\]
The equations of the asymptotes of the hyperbola bounding $\mathcal{X}_\delta$ are given by
$$
\Im(w) = \pm \frac{\sqrt{1-\delta^2}}{\delta}\Re(w).
$$
These asymptotes divide $\CP$ into four sectors with $\mathcal{X}_\delta$ belonging to the interior of the sector
\[\Big\{w\in\mathbb{C}: |\Im(w)| \leq \frac{\sqrt{1-\delta^2}}{\delta}\Re(w)\Big\}.\]
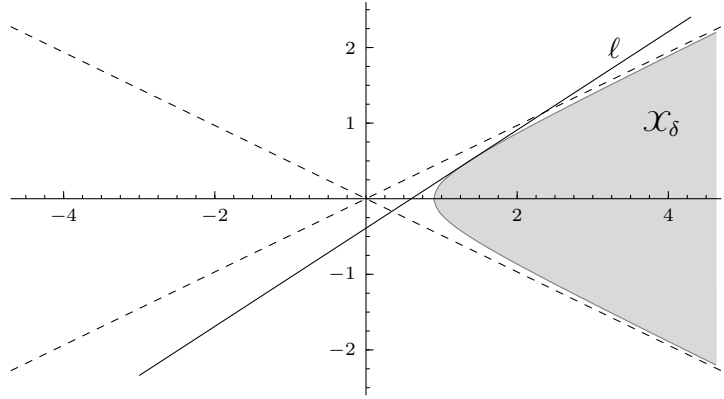
\begin{figure}[h]
\begin{center}
\begin{tikzpicture}
    \filldraw[gray!30] plot[domain=-2.2:2.2,smooth] ({0.9*sqrt(1+(\x)*(\x)/0.19)},{\x}) -- ({0.9*sqrt(1+(2.2*2.2/0.19)},-2.2);
    \draw[gray] plot[domain=-2.2:2.2,smooth] ({0.9*sqrt(1+(\x)*(\x)/0.19)},{\x});
    \draw (4.3,1) node[left]{$\mathcal{X}_\delta$};
    \draw (-4.7,0) -- (4.7,0);
    \foreach \a in {-4,-2,2,4}{
        \draw ({\a},0.1) -- ({\a},0) node[below, font=\tiny]{${\a}$};
    }
    \foreach \a in {-4.5,-4.25,...,4.5}{
        \draw ({\a},0.05) -- ({\a},0);
    }
    \draw (0,-2.6) -- (0,2.6);
    \foreach \a in {-2,-1,1,2}{
        \draw (0.1,{\a}) -- (0,{\a}) node[left, font=\tiny]{${\a}$};
    }
    \foreach \a in {-2.5,-2.25,...,2.5}{
        \draw ((0.05,{\a}) -- (0,{\a});
    }
    \draw[domain=-4.7:4.7, smooth, dashed] plot ({\x},{sqrt(.19)*(\x)/0.9});
    \draw[domain=-4.7:4.7, smooth, dashed] plot ({\x},{-sqrt(.19)*(\x)/0.9}); 
    \draw[domain=-3:4.3, smooth] plot ({\x},{sqrt((0.19)/(1.25))*(1.5*(\x)/0.9-1)});
    \draw (3.5,2) node[left]{$\mathcal{\ell}$};
\end{tikzpicture}
\caption{\small The region $\mathcal{X}_\delta$ corresponding to $\delta =\frac{9}{10}$, its asymptotes (dashed), and a tangent $\ell$ to its boundary.}\label{HyperbolaFig}
\end{center}
\end{figure}

We claim that the image of $\mathcal{X}_\delta$ under the inversion $\mathcal{I}(w) = \frac{1}{w-w_0}$ is a convex set if and only if $w_0$ lies in the closed sector opposite the one containing $\mathcal{X}_\delta$, which means
\begin{equation}\label{SOX}
|\Im(w_0)| \leq -\frac{\sqrt{1-\delta^2}}{\delta}\Re(w_0).
\end{equation}
 
 Suppose that $w_0\in \CP$ satisfies the preceding inequality.  Let $\mathcal{h}$ be the  hyperbola branch bounding $\mathcal{X}_\delta$,  let $\ell$ be a line tangent to $\mathcal{h}$, and let $\ell_{\mathcal{X}_\delta}\subset \CP$ be the closed halfplane  bounded by $\ell$ that contains all points of $\mathcal{X}_\delta$.   Observe that because $w_0$ lies in the closed sector opposite the one containing $\mathcal{X}_\delta$, $w_0$ does not belong to  $\ell_{\mathcal{X}_\delta}$.      Hence, the image of $\ell_{\mathcal{X}_\delta}$ under the linear-fractional inversion map $\mathcal{I}$ has the form $\mathcal{D}\setminus\{0\}$, where $\mathcal{D}$ is a closed disk having $0$ as a boundary point and where $\mathcal{I}(\mathcal{X}_\delta)\subseteq \mathcal{D}\setminus\{0\}$.   Because $\ell$ is an arbitrary tangent to $\mathcal{h}$, we conclude that $\mathcal{I}(\mathcal{X}_\delta)$ is the intersection of sets of the form $\mathcal{D}\setminus\{0\}$, where $\mathcal{D}$ is a closed disk having $0$ as a boundary point. Because an intersection of convex sets is convex and any closed disk having one or more boundary points removed is convex, we conclude $\mathcal{X}_\delta$ is convex.  
 
  Now suppose that $w_0\in \CP\setminus\mathcal{X}_\delta$ does not belong to the closed sector opposite the one containing $\mathcal{X}_\delta$. Then either (i) $w_0$ lies in the closed left halfplane $\{z: \Re(z) \le 0\}$ strictly above the negatively sloped asymptote $\mathcal{A}_{-}$ of $\mathcal{h}$ or $w_0$ lies in the open right halfplane $\{z: \Re(z) > 0\}$ and $\Im(w_0)\ge 0$, or (ii) $w_0$ lies in the closed left halfplane strictly below the positively sloped asymptote  $\mathcal{A}_{+}$ of $\mathcal{h}$  or $w_0$ lies in the open right halfplane and $\Im(w_0)< 0$. If (i) holds, let $L$ be the line through $w_0$ parallel to $\mathcal{A}_{-}$. If (ii) holds, let $L$ be the line through $w_0$ parallel to  $\mathcal{A}_{+}$. Observe that $L$ intersects $\mathcal{h}$ at exactly one point and that a slight rotation of $L$ about $w_0$, in the clockwise direction in case (i) and in the counter-clockwise direction in case (ii), yields a line $\tilde{L}$ through $w_0$ that intersects $\mathcal{h}$ at two points.  
  
  Let $P$ be the point in $\tilde{L}\cap\mathcal{h}$ that is farther from the origin. As illustrated in Figure~\ref{RayR} below, we let $\mathcal{R} \subset \CP$ be the ray along $\tilde{L}$ emanating from $P$ such that $\mathcal{R}$, except for $P$, lies entirely in the complement of $\mathcal{X}_\delta$. Observe that because $\tilde{L}$ passes through $w_0$, $\mathcal{I}(\mathcal{R})$ will be a half-closed line segment    containing its endpoint $\mathcal{I}(P)$ but not containing its other endpoint ${0}$.  Moreover, this line segment $\mathcal{I}(\mathcal{R})$ is contained in the complement of $\mathcal{I}(\mathcal{X}_\delta)$, except for its endpoint $\mathcal{I}(P)$.  Suppose, in order to obtain a contradiction that  $\mathcal{I}(\mathcal{X}_\delta)$ is convex. Then its closure is also convex. However, we claim the closure of $\mathcal{I}(\mathcal{X}_\delta)$ is simply  $\mathcal{I}(\mathcal{X}_\delta)\cup \{0\}$. Thus, accepting the validity of this claim, we see that the convexity of   $\mathcal{I}(\mathcal{X}_\delta)$ implies the convexity of $\mathcal{I}(\mathcal{X}_\delta) \cup \{0\}$.  However, the line segment connecting $\mathcal{I}(P) \in  (\mathcal{I}(\mathcal{X}_\delta) \cup \{0\})$ to $0\in  (\mathcal{I}(\mathcal{X}_\delta) \cup \{0\})$ lies, except for its endpoints in the complement of  $\mathcal{I}(\mathcal{X}_\delta) \cup \{0\}$, a contradiction.

 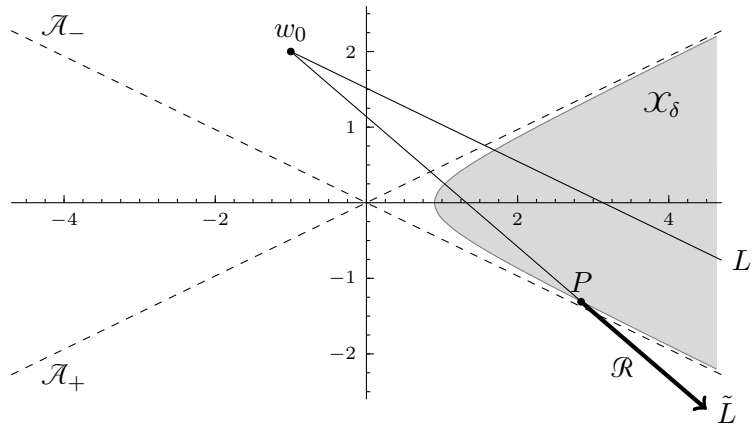
\begin{figure}[h]
\begin{center}
\begin{tikzpicture}
    \filldraw[gray!30] plot[domain=-2.2:2.2,smooth] ({0.9*sqrt(1+(\x)*(\x)/0.19)},{\x}) -- ({0.9*sqrt(1+(2.2*2.2/0.19)},-2.2);

    \draw[gray] plot[domain=-2.2:2.2,smooth] ({0.9*sqrt(1+(\x)*(\x)/0.19)},{\x});
    \draw (4.3,1.3) node[left]{$\mathcal{X}_\delta$};
    \draw (-4.7,0) -- (4.7,0);
    \foreach \a in {-4,-2,2,4}{
        \draw ({\a},0.1) -- ({\a},0) node[below, font=\tiny]{${\a}$};
    }
    \foreach \a in {-4.5,-4.25,...,4.5}{
        \draw ({\a},0.05) -- ({\a},0);
    }
    \draw (0,-2.6) -- (0,2.6);
    \foreach \a in {-2,-1,1,2}{
        \draw (0.1,{\a}) -- (0,{\a}) node[left, font=\tiny]{${\a}$};
    }
    \foreach \a in {-2.5,-2.25,...,2.5}{
        \draw ((0.05,{\a}) -- (0,{\a});
    }

    \draw[domain=-4.7:4.7, smooth, dashed] plot ({\x},{sqrt(.19)*(\x)/0.9});
    \draw (-4,-2) node[below]{$\mathcal{A}_{+}$};

    \draw[domain=-4.7:4.7, smooth, dashed] plot ({\x},{-sqrt(.19)*(\x)/0.9});
    \draw (-4,2) node[above]{$\mathcal{A}_{-}$}; 

    \draw[domain=-1:4.7, smooth] plot ({\x},{-sqrt(.19)*(\x+1)/0.9 + 2}) node[right]{$L$}; %%% The line L

    \filldraw[black] (-1,2) circle[radius=1.25pt];
    \draw (-1,2) node[above]{$w_0$}; %%% w_0 = -1 + 2i

    \filldraw[black] (2.84,-1.31) circle[radius=1.25pt];
    \draw (2.84,-1.31) node[above]{$P$}; %%% P = 2.84 - 1.31i

    \draw[domain=-1:4.5, smooth] plot ({\x},{-sqrt(0.60)*(\x+1)/0.9 + 2}) node[right]{$\tilde{L}$}; %%% The line \tilde{L}

\draw[-to, ultra thick] (2.84,-1.31) -- (4.5,-2.7336); %%% The ray \mathcal{R}

\draw (3.4,-1.9) node[below]{$\mathcal{R}$};

\end{tikzpicture}
\caption{\small The lines $L$ and $\tilde{L}$, as well as the ray $\mathcal{R}$.}\label{RayR}
\end{center}
\end{figure}

  We complete the proof by verifying the claim made in the preceding paragraph that the closure of $\mathcal{I}(\mathcal{X}_\delta)$ is  $\mathcal{I}(\mathcal{X}_\delta)\cup \{0\}$.  That $0$ belongs to the closure of $\mathcal{I}(\mathcal{X}_\delta)$ is clear because $\mathcal{X}_\delta$ is unbounded, and, for any unbounded sequence $(z_n)$  in $\mathcal{X}_\delta$, $0$ is limit point of $\{\mathcal{I}(z_n):n\in \mathbb{N}\}$. To complete the verification of our claim, we show that any nonzero limit point of $\mathcal{I}(\mathcal{X}_\delta)$ belongs to $\mathcal{I}(\mathcal{X}_\delta)$. Let $z_*$ be a nonzero limit point of $\mathcal{I}(\mathcal{X}_\delta)$, so that there is a sequence $(z_n)$ in $\mathcal{X}_\delta$ such that $\lim (\mathcal{I}(z_n)) =z_*$. Because $z_*\ne 0$, the sequence $(z_n)$ is bounded.  Thus $(z_n)$ has a convergence subsequence $(z_{j_n})$ whose limit $z$ must belong to  $\mathcal{X}_\delta$ because $\mathcal{X}_\delta$ is closed.  Because $\mathcal{I}$ is continuous on $\mathcal{X}_\delta$ ($w_0\in \CP\setminus\mathcal{X}_\delta$),  $z_* = \lim(\mathcal{I}(z_{j_n})) =\mathcal{I}(z)\in \mathcal{I}(\mathcal{X}_\delta)$, as desired.  \end{proof}

 Let $w_1$ and $w_2$ be distinct points in $\CP$ and let $0< d < |w_1 - w_2|$.  Observe that 
\begin{equation}\label{GHB}
 \mathcal{h}: =\{w\in \CP: |w -w_1| - |w - w_2| = d \}
\end{equation}
is one branch of the hyperbola with foci $w_1$ and $w_2$ and constant difference $d$; specifically, the branch having vertex 
$$
\frac{w_1 + w_2}{2}   + \frac{d}{2}\frac{w_2 - w_1}{|w_2-w_1|}.
$$
Note that the complement of  $\mathcal{h}$  has two open connected components, one contains the focus $w_2$ and the other contains the center $(w_1 + w_2)/2$.  

 \begin{theorem}  Let $\ph(z)=az+b$ be a selfmap of $\D$ and let $d = \sqrt{|1-a|^2-|b|^2}$.
    Then $B(C_\phi)$ is convex if and only if
    \begin{equation}\label{CFC}
     2d|\Im(a)| \le |b|(1-|a|^2).
     \end{equation}
     \end{theorem}
    \begin{proof}
We already noted (in the paragraph preceding Lemma~\ref{L:hyperbola_inversion}) that if $b=0$, then the theorem holds. Thus, to complete the proof, we consider the case $b\ne 0$. There are two possibilities $a + |b| =1$ or $a + |b| \ne 1$. In the former case, (\ref{CFC}) is satisfied (because $\Im(a) = 0$) and $B(C_\phi)$ is convex by part (iii) of Theorem~\ref{T:Berezin_range_az_b} because $a + |b| =1$ is equivalent to the fixed point $b/(1-a)$ of $\phi$ belonging to $\partial \D$ (Proposition~\ref{LSMBFP}). Thus, to complete the proof, we need only consider the case for which $b\ne 0$ and $\phi$ fixes a point inside $\D$; equivalently, $a + |b| \ne 1$.

  Let $b\ne 0$ and assume $a+|b| \ne 1$.   By Theorem~\ref{T:Berezin_range_az_b}, 
  $$
  B(C_{\phi}) = \left\{z\in \CP\setminus\{0\}: \left|\frac{1}{z} - a\right| - \left|\frac{1}{z} -1\right| \ge d\right\},
  $$
so that the boundary of the set of reciprocals of points in  $B(C_\phi)$  is the hyperbola branch $\mathcal{h}$ of (\ref{GHB}) where $w_1=a$ and $w_2= 1$ are the hyperbola's foci. The hyperbola is not degenerate because $d = \sqrt{|1-a|^2-|b|^2} < |1-a|$ ($b\ne0$).  Thus, $B(C_\phi)$ itself consists of reciprocals of points either on the hyperbola branch $\mathcal{h}$ or in the component of its complement containing the focus 1. Let $\tilde{\mathcal{h}}$ be the union of $\mathcal{h}$ and the component of its complement containing the focus 1. 
 
  Via the transformation of $\CP$ defined by
 $$
 T(z) = \left(z-\frac{1+a}{2}\right)\frac{|1-a|}{1-a},
 $$
  we translate and rotate $\mathcal{h}$ into standard position, so that it becomes the right-hand branch of the hyperbola with center  $0$, foci $\pm |1-a|/2$, and constant difference $d$.   Hence,
  \begin{align*}
  T(\tilde{\mathcal{h}}) & = \left\{w\in \CP: \left|\frac{2}{|1-a|}w + 1 \right| - \left|\frac{2}{|1-a|}w -1\right| \ge  \frac{2d}{|1-a|}\right\}\\
  & = \left\{w\in \CP: \frac{2}{|1-a|}w\in \mathcal{X}_{d/|1-a|}\right\},
  \end{align*}
  where $\mathcal{X}_{d/|1-a|}$ is defined by Equation~(\ref{HBX}) of Lemma~\ref{L:hyperbola_inversion}. 
  It follows that 
  $$
  B(C_\phi) = \left\{z\in\CP\backslash\{0\}: \frac{1}{z}\in \tilde{\mathcal{h}}\right\} = \left\{z\in\CP\backslash\{0\}: \frac{2}{|1-a|}T(1/z)\in \mathcal{X}_{d/|1-a|}\right\}.
$$
   A little algebra shows that the preceding equation is equivalent to 
  \begin{equation}\label{BRC}
   B(C_\phi) = \left\{\frac{2}{1-a}\frac{1}{w+ \frac{1+a}{1-a}}: w\in  \mathcal{X}_{d/|1-a|}\right\}.
\end{equation}
   Equation \eqref{BRC} shows that $B(C_\phi)$ is convex if only if the image of $\mathcal{X}_{d/|1-a|}$ under $
   w\mapsto \frac{1}{w+ \frac{1+a}{1-a}}$ is convex (because for $c\in \CP\setminus\{0\}$  and $S\subseteq \CP$, $cS$ is convex if and only if $S$ is convex). By Lemma~\ref{L:hyperbola_inversion}, with $\delta =d/|1-a|$ and $w_0 = - \frac{1+a}{1-a}$,   we conclude that  $B(C_\phi)$ is convex if and only if  
 \begin{equation}\label{CCEQ}
 \Big|\Im\Big(-\frac{1+a}{1-a}\Big)\Big| + \frac{\sqrt{1-\delta^2}}{\delta}\Re\Big(-\frac{1+a}{1-a}\Big) \leq 0.
 \end{equation}
Because
\[\frac{\sqrt{1-\delta^2}}{\delta} = \sqrt{\frac{1}{\delta^2}-1} = \sqrt{\frac{|1-a|^2}{d^2}-1} = \frac{|b|}{d},\]
and
\[\frac{1+a}{1-a} = \frac{(1+a)(1-\overline{a})}{|1-a|^2} = \frac{1-|a|^2}{|1-a|^2}+\frac{2i\Im(a)}{|1-a|^2},\]
we obtain from inequality~(\ref{CCEQ}) that
\[\frac{2|\Im(a)|}{|1-a|^2}-\frac{|b|}{d}\cdot\frac{1-|a|^2}{|1-a|^2} \leq 0,\]
which is equivalent to $2d|\Im(a)| \le |b|(1-|a|^2)$, as desired.
\end{proof}

We conclude with a couple of natural questions raised by our work:
\begin{enumerate}
\item Proposition~\ref{symprop} shows that a sufficient condition for $B(C_\phi)$ to be symmetric relative to the real axis is  that $\phi$ map a diameter of the $\D$ into itself.  Is this condition necessary? 
\item  For an arbitrary analytic selfmap $\phi$ of $\D$, are $0$ and $1$ the only possible corner points of $B(C_\phi)$?
\end{enumerate}

\section*{Acknowledgments}

Our initial equations and inequalities for the $b\ne 0$ cases of Theorem 4.1 were expressed in terms of  $x= \Re(z)$ and $y=\Im(z)$. In this context, Claude Opus 5.5 provided a factorization corresponding to that on line 7 of Case 2 of the proof of Theorem 4.1. ChatGPT-5 produced tikz code yielding Figures 1 and 2.  In addition, the authors used AI-assisted searches to find results and references relevant to our work.  The preceding represents the full extent of the authors' AI use in preparing this paper, for which the authors take full responsibility. 

The authors used Mathematica 9.0.1 (to, e.g., produce Figures 3-6).

Trieu Le is partially supported by Simons Foundation Travel Support for Mathematicians MPS-TSM-00002303.

{\footnotesize{

    \bibliography{citations}
    \bibliographystyle{plain}
    }}

\end{document}